\documentclass{amsart}
\usepackage[english]{babel}
\usepackage[latin1]{inputenc}

\usepackage{amsmath}
\usepackage{amsfonts}
\usepackage{amssymb}
\usepackage{amsthm}
\usepackage{graphicx,epsfig}
\usepackage{amscd}
\DeclareMathOperator{\re}{Re}
\DeclareMathOperator{\im}{Im}
\DeclareMathOperator{\Res}{Res}

\newtheorem{defn}{Definition}
\newtheorem{thm}{Theorem}
\newtheorem{cor}{Corollary}
\newtheorem{prop}{Proposition}

\newtheorem{lemma}{Lemma}
\newtheorem{rem}{Remark}
\newtheorem{ex}{Example}

\newcommand{\R}{\mathbb{R}}
\newcommand{\h}{\mathbb{H}}
\newcommand{\Z}{\mathbb{Z}}
\newcommand{\N}{\mathbb{N}}
\newcommand{\C}{\mathbb{C}}

\newcommand{\rmd}{\mathrm{d}}
\newcommand{\rmO}{\mathrm{O}}
\newcommand{\rmo}{\mathrm{o}}

\newcommand{\cA}{{\mathcal A}}
\newcommand{\cB}{{\mathcal B}}
\newcommand{\cC}{{\mathcal C}}
\newcommand{\cD}{{\mathcal D}}
\newcommand{\cZ}{{\mathcal Z}}

\begin{document}

\title[Flux for Bryant Surfaces]
{Flux for Bryant Surfaces
and Applications to Embedded
  Ends of Finite Total Curvature}
\author{Beno\^\i t Daniel}
\date{}

%\begin{abstract}
%We give a geometric description of embedded ends of constant mean
%curvature $1$ and of finite total curvature in hyperbolic space. In
%particular, we show 
%that we can define an axis for these ends that are asymptotic to a
%catenoid cousin. We also compute the flux of Killing fields through
%these ends, and we deduce some geometric properties and some
%analogies with minimal surfaces in Euclidean
%space. Finally, we answer partially a conjecture of Rossman, Umehara 
%and Yamada by showing that the flux matrix they have defined is
%equivalent to the flux of Killing fields in the case of embedded ends
%of finite total curvature.
%\end{abstract}

\begin{abstract}
We compute the flux of Killing fields through ends of constant mean curvature $1$ in hyperbolic space, and we prove a result conjectured by Rossman, Umehara and Yamada : the flux matrix they have defined is equivalent to the flux of Killing fields. We next give a geometric description of embedded ends of finite total curvature. In particular, we show 
that we can define an axis for these ends that are asymptotic to a
catenoid cousin. We also compute the flux of Killing fields through
these ends, and we deduce some geometric properties and some
analogies with minimal surfaces in Euclidean space.
\end{abstract}


\maketitle
%\tableofcontents

\section{Introduction}

Bryant surfaces are surfaces with constant mean curvature one in 
hyperbolic 3-space $\h^3$ (with the convention that the mean curvature
of a surface is one half of the trace of its second fundamental form). These surfaces have been studied first by
Bryant (\cite{bryant}). He derived a representation in terms of
holomorphic data, analogous to the Weierstrass data for minimal
surfaces in $\R^3$.

Umehara and Yamada have defined the notion of regular ends of Bryant
surfaces (\cite{umehara}): these are ends conformally parametrized by
the punctured 
complex disk and such that the hyperbolic Gauss map extends
meromorphically to the puncture (if the hyperbolic Gauss map has an
essential singularity at the puncture, the end is said to be
irregular). They also have studied the 
Weierstrass data of Bryant surface ends of finite total curvature.

\begin{figure}[htbp]
\begin{center}
\includegraphics{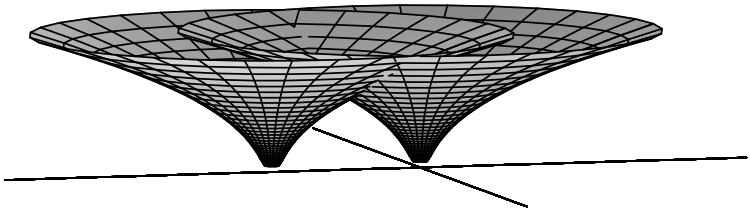}
\caption{Two half-catenoid cousins asymptotic in the sense of S\'a
  Earp and Toubiana but with different axes.} \label{twocatenoids}
\end{center}
\end{figure}

Collin, Hauswirth and Rosenberg (\cite{collin}) have shown that properly
embedded annular ends have finite total curvature and are regular. Yu
(\cite{yu}) has shown that irregular ends are never embedded.
S\'a Earp and Toubiana (\cite{toubiana}) studied the geometry of
embedded ends of finite total curvature (hence regular). They showed
that, in the upper half-space model of $\h^3$, such ends are, up to an
isometry of $\h^3$, vertical Euclidean graphs and are asymptotic to a
catenoid cousin of revolution or a horosphere as vertical Euclidean
graphs. They also defined the growth of such ends. If $E$ is a
half-catenoid cousin whose asymptotic boundary is $\infty$, then the
image of $E$ by a Euclidean horizontal translation
(which is a parabolic isometry of $\h^3$) is asymptotic to $E$ in the
sense of S\'a Earp and Toubiana (see figure \ref{twocatenoids}).

There exist two notions of flux for Bryant surfaces. The first flux is 
the flux of Killing fields. This flux was introduced
by Korevaar, Kusner, Meeks and Solomon (\cite{kkms}) as an analogue
of the flux defined by Korevaar, Kusner and Solomon for constant mean
curvature surfaces in $\R^3$ (\cite{kks}). It is the sum of an
integral along a curve $\Gamma$
% that is a generator of the fundamental
%group of the end, 
and of an integral over a compact surface whose boundary is $\Gamma$.
% (actually we will show in lemma \ref{integral} that we can
%reduce ourselves to one single integral along $\Gamma$). 
This flux is
a homology invariant. The second flux is the residue-type flux matrix defined by Rossman, Umehara and Yamada (\cite{umeharaflux}). This flux can be easily computed from the Bryant representation of the surface. It is also a homology invariant. Rossman, Umehara and Yamada conjectured that these two notions of flux were equivalent.

In this paper, we prove this conjecture. We compute the flux of Killing fields associated to translations and rotations through Bryant surface ends. We show that it only depends on the residues of three meromorphic one-forms (theorems \ref{fluxtranslation} and \ref{fluxrotation}). These residues are, up to constant factors, the coefficients of the flux matrix defined by Rossman, Umehara and Yamada. Moreover, we define a complex polynomial of degree at most two, called flux polynomial, whose coefficients are these residues (theorem \ref{fluxpolynomial}). This polynomial contains all the information given by the
flux and satisfying a ``balancing formula''.

The second aim of this paper is to complete the geometric study of embedded ends of finite total curvature started in \cite{toubiana}. We show that we can 
define an axis for such ends that are asymptotic to a catenoid
cousin (theorem \ref{uniqueness}). This means that these ends are
asymptotically surfaces of 
revolution. We call these ends catenoidal ends. The analogous result
for embedded ends of finite total curvature of minimal surfaces in $\R^3$
has been proved by Schoen (\cite{schoen}).

We next compute the flux for embedded ends of finite total curvature.
We obtain that the flux of the Killing field associated to the
translation along the geodesic $(\cC,\cD)$ through a catenoidal end is
$$\varphi=\pi(1-\mu^2)(2\re(\cA,\cC,\cD,\cB)-1)$$ where $(\cA,\cB)$ is
the axis of the end, $1-\mu$ its growth, and where $(\cA,\cC,\cD,\cB)$
denotes the cross-ratio (theorem \ref{theoremfluxcatenoidal}). This
formula is one of the simplest we could expect, since it depends only
on the asymptotic behaviour of the end.
We also show that the flux for a horospherical end is zero if and only
if its Hopf differential is regular at the end
(theorem \ref{theoremfluxhorospherical}).

Thus, the flux for Bryant surfaces plays the same role as the flux and the
torque for minimal surfaces in $\R^3$ (the torque is defined in
\cite{kk} ; see also \cite{hoffman}
for definitions and basic properties of the flux and the
torque). Indeed, the flux 
and the torque for a catenoidal end depend only on the growth and the
axis of the end, and the torque for a planar end (the analogue of a
horospherical end) is zero if and only if the Hopf differential is
regular at the end, {\it i.e.} the degree of the Gauss map
at the end is at least $3$ (see \cite{romon}).

%We show that we can associate to an end a complex polynomial of degree
%at most two containing all the information given by the
%flux and satisfying a ``balancing formula'' (theorem \ref{thmpoly} and
%corollary \ref{poly}). 
Finally, we give some geometric applications of the flux.
If a Bryant surface has exactly two catenoidal ends (and
no others) with distinct asymptotic boundaries, then the ends have the
same growth and the same axis (proposition \ref{deux}). If a Bryant
surface has exactly three catenoidal ends (and no others) with
distinct asymptotic boundaries, then the axes are coplanar and
concurrent (possibly in the asymptotic boundary of $\h^3$)
(proposition \ref{trois}). The same results hold for minimal surfaces
in $\R^3$.

%Finally we compute the residue-type flux matrix defined
%by Rossman, Umehara and Yamada in \cite{umeharaflux} for embedded
%ends of finite total curvature. We obtain that the flux matrix
%and the flux polynomial are equivalent for these ends. Thus we answer
%positively a conjecture of Rossman, Umehara and Yamada in this
%special case.

\section{Preliminaries and notations} \label{prel}

In this paper, the model used for hyperbolic 3-space is the upper
half-space model :
$$\h^3=\{(u,v,w)\in\R^3|w>0\}=\{(\zeta,w)\in\C\times\R|w>0\}$$ 
with the metric 
$$\rmd s^2=\frac{\rmd u^2+\rmd v^2+\rmd w^2}{w^2}
=\frac{|\rmd\zeta|^2+\rmd w^2}{w^2}.$$
$<,>$ and $||.||$ denote respectively the hyperbolic metric and the
hyperbolic norm on $\h^3$. If $X_1=(\alpha_1,\beta_1)$ and
$X_2=(\alpha_2,\beta_2)$ are two vectors in the tangent space of
$\h^3$ at the point $(\zeta,w)$, then 
$<X_1,X_2>=(\re(\bar{\alpha}_1\alpha_2)+\beta_1\beta_2)/w^2$. 

In the model of the unit ball of $\R^3$ for hyperbolic space, the
asymptotic boundary of hyperbolic space is the sphere of radius
$1$. In the half-space model, we identify the asymptotic boundary of
$\h^3$ with the Riemann sphere $\bar{\C}$ composed of the plane
$\{\zeta=0\}$ and of the point at infinity which we denote $\infty$.

The asymptotic boundary of a part of $\h^3$ is the set of its
accumulation points in $\bar{\C}$.

The identification between the upper half-space model and 
the Minkowski model for $\h^3$ is the same as that described in
\cite{toubiana} (remark 1.11). Consequently, if $f$ is a constant mean curvature one immersion of a Riemann surface $M$ into the Minkowski model of the hyperbolic space, if $F=\left( 
\begin{array}{cc}
A & B \\
C & D
\end{array} \right)$ is its Bryant representation (see
\cite{bryant}), then we have $f=FF^*$, and the corresponding immersion $X=(\zeta,w):M\to\h^3$ in the
upper half-space model is given by 
\begin{equation} \label{zetaabcd}
\zeta=\frac{\bar{A}C+\bar{B}D}{|A|^2+|B|^2}
\end{equation}
and
\begin{equation} \label{wabcd}
w=\frac{1}{|A|^2+|B|^2}.
\end{equation}
We recall that $A$, $B$, $C$ and $D$ are holomorphic functions defined on the universal cover of $M$ and satisfying $AD-BC=1$ and $\rmd A\rmd D-\rmd B\rmd C=0$.

If $(g,\omega)$ denote the Weierstrass data of the end (see \cite{bryant} or \cite{umehara}), the $2$-form
$\omega\rmd g$ is called the Hopf differential of the end. It is single-valued on $M$ (contrarily to $g$ and $\omega$). It is
invariant by an isometry of $\h^3$. 

The hyperbolic Gauss map is given by $G=\frac{\rmd C}{\rmd A}=\frac{\rmd D}{\rmd B}$. It is single-valued on $M$. This expression slightly differs from that of \cite{bryant}, \cite{umehara} and other papers because of the chosen identification (see \cite{toubiana}, remark 1.11). The one-form $\omega^\#=-\frac{\omega\rmd g}{\rmd(1/G)}$ is also single-valued on $M$ (the couple $(1/G,\omega^\#)$ gives the Weierstrass data of the dual immersion, see \cite{umeharadual}). Hence the following one-forms are single-valued on $M$ :
$$B\rmd A-A\rmd B=-\frac{\omega^\#}{G^2}=\frac{\omega\rmd g}{\rmd G},$$
$$C\rmd B-B\rmd A=\frac{\omega^\#}{G}=-G\frac{\omega\rmd g}{\rmd G},$$
$$D\rmd C-C\rmd D=-\omega^\#=G^2\frac{\omega\rmd g}{\rmd G}.$$

For regular ends of finite total
curvature, the Hopf differential $\omega\rmd g$ has a pole of order greater than or equal to $-2$ at
zero (see \cite{umehara}). Its order does not depend on the parametrization. 

In the Minkowski model, a direct isometry is a map $N\mapsto PNP^*$ where
$P=\left( 
\begin{array}{cc}
\alpha & \beta \\
\gamma & \delta
\end{array} \right)\in{\mathrm{SL}}_2(\C)$.
In the half-space model, this isometry induces on $\bar{\C}$ the map 
$\zeta\mapsto\frac{\delta\zeta+\gamma}{\beta\zeta+\alpha}$ because of
the chosen identification. 

If ${\mathcal A}$ and ${\mathcal B}$ are two distinct points in
$\bar{\C}$, 
$({\mathcal A},{\mathcal B})$ denotes the oriented geodesic of $\h^3$ going
from ${\mathcal A}$ to ${\mathcal B}$.

If $z_1$, $z_2$, $z_3$ and $z_4$ are four points in $\bar{\C}$ such
that $z_1\neq z_4$ and $z_2\neq z_3$, we define their cross-ratio by
$$(z_1,z_2,z_3,z_4)=\frac{z_3-z_1}{z_3-z_2}\cdot\frac{z_4-z_2}{z_4-z_1}.$$
We recall that there exists a direct isometry (respectively an
indirect isometry) of $\h^3$ which maps
$z_1$, $z_2$, $z_3$ and $z_4$ to $z'_1$, $z'_2$, $z'_3$ and $z'_4$
respectively (where $z'_1\neq z'_4$ and $z'_2\neq z'_3$) if and only
if $(z_1,z_2,z_3,z_4)=(z'_1,z'_2,z'_3,z'_4)$ (respectively
$(z_1,z_2,z_3,z_4)=\overline{(z'_1,z'_2,z'_3,z'_4)}$). 

In this paper, $\Omega$ will denote any neighbourhood of $0$ in $\C$,
and $\Omega^*$ will denote the set $\Omega\setminus\{0\}$. 

The flux of a Killing field $Y$ through an annular Bryant surface end
$E$ is defined by 
$$\varphi=\int_\Gamma<\eta,Y>-2\int_K<\nu,Y>$$ where $\Gamma$ is a
generator of $\pi_1(E)$, $K$ a topological disk whose boundary is
$\Gamma$, $\eta$ the conormal to $\Gamma$ in the direction of the
asymptotic boundary of the end and $\nu$ the normal to $K$ chosen as
follows : we choose on $\Gamma$ the orientation such that
$(\Gamma,\eta,-\vec{H})$ is the orientation of $\h^3$ and $\nu$ such
that it induces the same orientation on $\Gamma$.
The normal $\nu$ induces an orientation on $K$ and
$\Gamma$. These choices have been made in order to be compatible with
Stokes's forumula.

This number $\varphi$ does not depend on the choices of $\Gamma$
and $K$ (see \cite{kks} and \cite{kkms}). We shall notice that in
\cite{kks} and \cite{kkms} the 
mean curvature is defined as the trace of the second fundamental form
(and not its half), which explains the coefficient $2$ in the formula.

If $\alpha$ is a $n$-form on $\h^3$ and $X$ a vector field, then the
interior product of $\alpha$ by $X$ is denoted ${\rm i}_X\alpha$ and
defined by 
${\rm i}_X\alpha(\xi_1,\dots,\xi_{n-1})
=\alpha(X,\xi_1,\dots,\xi_{n-1})$. The Lie derivative of $\alpha$ with
respect to $X$ is denoted ${\rm L}_X\alpha$. We recall Cartan's
forumula : ${\rm L}_X\alpha=\rmd({\rm i}_X\alpha)+{\rm i}_X\rmd\alpha$.

%Finally, $[\rho^\alpha]$ denotes any function of $\rho$ that is
%equivalent to $\beta\rho^\alpha$ with $\beta\neq 0$ when $\rho$ tends to zero.

\section{Flux of Killing fields}

\subsection{Killing fields associated to translations}

\begin{defn}
Let ${\mathcal A}$ and ${\mathcal B}$ be two distinct points in
$\bar{\C}$. Let 
$\Phi_t$ be the translation of distance $t\in\R$ along the geodesic
$({\mathcal A},{\mathcal B})$. Then the vector field $Y$ defined by
$$\frac{\rmd\Phi_t}{\rmd t}=Y(\Phi_t)$$ is called the Killing field
associated to the translation along $({\mathcal A},{\mathcal B})$.
\end{defn}

The Killing field associated to the translation along $(\cB,\cA)$ is the opposite of the Killing field associated to the translation along $(\cA,\cB)$.
Elementary computations give the following lemma.

\begin{lemma}
The Killing field associated to the translation along $(0,\infty)$ is
$$Y(\zeta,w)=(\zeta,w).$$
\end{lemma}

\begin{lemma} \label{killingtrans}
Let $\zeta_0\in\C$. The Killing field associated to the
translation along $(\zeta_0,0)$ is
\begin{displaymath}
Y(\zeta,w)=\left(
\begin{array}{c}
-\frac{w^2}{\bar{\zeta}_0}+\frac{\zeta^2}{\zeta_0}-\zeta \\
2w\re\frac{\zeta}{\zeta_0}-w
\end{array}
\right).
\end{displaymath}
\end{lemma}

\begin{proof}
The map $$\Phi:(u',v',w')\mapsto\frac{1}{{u'}^2+{v'}^2+{w'}^2}(u',v',w')$$
is an isometry of $\h^3$ which maps the geodesic $(\zeta'_0,\infty)$
(where $\zeta'_0=\zeta_0/|\zeta_0|^2$) to the geodesic $(\zeta_0,0)$. Hence the
Killing field associated to the translation along $(\zeta_0,0)$ is
given by $Y(P)=\Phi_*Z(P)=\rmd_{\Phi^{-1}(P)}\Phi\cdot Y(\Phi^{-1}(P))$ for
each $P=(u,v,w)=(\zeta,w)\in\h^3$, where $Z$ is the
Killing field associated to the translation along $(\zeta'_0,\infty)$.

We have $\Phi^{-1}(P)=(u',v',w')$ where $u=ru'$, $v=rv'$, $w=rw'$ and
$r=u^2+v^2+w^2$. Hence we have $Z(\Phi^{-1}(P))=(u'-u'_0,v'-v'_0,w')$
with $\zeta'_0=u'_0+iv'_0$, so
\begin{eqnarray*}
Y(P) & = & \frac{1}{({u'}^2+{v'}^2+{w'}^2)^2}\\
& & \times\left(
\begin{array}{ccc}
{v'}^2+{w'}^2-{u'}^2 & -2u'v' & -2u'w' \\
-2u'v' & {u'}^2+{w'}^2-{v'}^2 & -2v'w' \\
-2u'w & -2v'w' & {u'}^2+{v'}^2-{w'}^2
\end{array} \right) \\
& & \times\left(
\begin{array}{c}
u'-u'_0 \\
v'-v'_0 \\
w'\end{array} \right) \\
& = & r^2\left(
\begin{array}{c}
-u'_0({v'}^2+{w'}^2-{u'}^2)+2v'_0u'v'-u'{v'}^2-u'{w'}^2-{u'}^3 \\
2u'_0u'v'-v'_0({u'}^2+{w'}^2-{v'}^2)-{u'}^2v'-v'{w'}^2-{v'}^3 \\
2u'_0u'w'+2v'_0v'w'-{u'}^2w'-{v'}^2w'-{w'}^3
\end{array} \right) \\
& = & \left(
\begin{array}{c}
u'_0(u^2-v^2-w^2)+2v'_0uv-u \\
2u'_0uv+v'_0(v^2-u^2-w^2)-v\\
2u'_0uw+2v'_0vw-w
\end{array}
\right)
=\left(
\begin{array}{c}
-\frac{w^2}{\bar{\zeta}_0}+\frac{\zeta^2}{\zeta_0}-\zeta \\
2w\re\frac{\zeta}{\zeta_0}-w
\end{array}
\right).
\end{eqnarray*}
\end{proof} 

\subsection{Killing fields associated to rotations}

\begin{defn}
Let $\cA$ and $\cB$ be two distinct points in $\bar{\C}$. Let
$R_\theta$ be the rotation of angle $\theta$ (in the direct sense)
about the geodesic $(\cA,\cB)$. Then the vector field $Y$ defined by
$$\frac{\rmd R_\theta}{\rmd\theta}=Y(R_\theta)$$ is called the Killing field
associated to the rotation about $(\cA,\cB)$.
\end{defn}

The Killing field associated to the rotation about $(\cB,\cA)$ is the opposite of the Killing field associated to the rotation about $(\cA,\cB)$.
Elementary computations give the following lemma.

\begin{lemma}
The Killing field associated to the rotation about $(0,\infty)$ is
$$Y(\zeta,w)=(i\zeta,0).$$
\end{lemma}

\begin{lemma} \label{killingrot}
Let $\zeta_0\in\C$. The Killing field associated to the
rotation about $(\zeta_0,0)$ is
\begin{displaymath}
Y(\zeta,w)=\left(
\begin{array}{c}
i\frac{w^2}{\bar{\zeta}_0}
+i\frac{\zeta^2}{\zeta_0}-i\zeta \\
-2w\im\frac{\zeta}{\zeta_0}
\end{array}
\right).
\end{displaymath}
\end{lemma}

\begin{proof}
We proceed as for lemma \ref{killingtrans} and we use the same
notations. Since the map $\Phi$ is an indirect isometry of $\h^3$, we
have $Y=-\Phi_*Z$ where $Z$ is the
Killing field associated to the rotation about $(\zeta'_0,\infty)$.
\end{proof}

\subsection{Flux of Killing fields associated to translations} \label{computationtranslation}

In this section, $\zeta_0$ and $\zeta_1$ are two complex numbers such that $\zeta_0\neq 0$, and $E$ denotes a Bryant surface end 
whose Bryant representation is $F=\left(
\begin{array}{cc}
A & B \\
C & D
\end{array} \right):\Omega^*\to\mathrm{SL}_2(\C)$. We denote by
$X=(\zeta,w):\Omega^*\to\h^3$ the corresponding conformal immersion in
the upper half-space model.
%and $Y$ is the Killing field associated to the translation along
%the geodesic $(\zeta_0+\zeta_1,\zeta_1)$.

We will denote by $(\rho,\tau)$ the polar coordinates in $\Omega$
({\it i.e.} $z=\rho e^{i\tau}$). We have the following
relationships for derivation operators :
$$\frac{\partial}{\partial z}=
\frac{1}{2z}\left(\rho\frac{\partial}{\partial\rho}
-i\frac{\partial}{\partial\tau}\right),$$
$$\frac{\partial}{\partial\bar{z}}=
\frac{1}{2\bar{z}}\left(\rho\frac{\partial}{\partial\rho}
+i\frac{\partial}{\partial\tau}\right).$$

\begin{lemma} \label{integral}
%Let $E$ be a Bryant surface end given 
%by a conformal immersion $X=(\zeta,w):\Omega^*\to\h^3$ in
%the upper half-space model. Let
%$\zeta_0$ and $\zeta_1$ be two complex numbers, with $\zeta_0\neq 0$, 
Let $Y$ be the Killing field associated to the translation along
the geodesic $(\zeta_0+\zeta_1,\zeta_1)$.
Then if $\rho$ is a sufficiently small positive number, the flux of
$Y$ through $E$ is 
$$\varphi=\int_0^{2\pi}\re(s(\rho,\tau))\rmd\tau$$
where
\begin{eqnarray*}
s(\rho,\tau) & = &
\frac{\zeta-\zeta_1}{w^2}
\overline{\left(\rho\frac{\partial\zeta}{\partial\rho}
+i\frac{\partial\zeta}{\partial\tau}\right)}
\left(1-\frac{\zeta-\zeta_1}{\zeta_0}\right) \\
& & +\frac{\rho}{\zeta_0}\frac{\partial\zeta}{\partial\rho}
-\frac{2i}{\zeta_0}\frac{\partial\zeta}{\partial\tau}\ln w
-\frac{\rho}{w}
\frac{\partial w}{\partial\rho}
\left(2\frac{\zeta-\zeta_1}{\zeta_0}-1\right).
\end{eqnarray*}
\end{lemma}

\begin{proof}
Let $\rho>0$ such that the circle $\{z\in\C||z|=\rho\}$ is included in
$\Omega$. Let $\Gamma$ be the curve on $E$ defined by $\tau\mapsto
X(\rho e^{i\tau})$. Let $K$ be a disk whose boundary is $\Gamma$.

We remark that we must take on $\Gamma$ the orientation given by
$-\Gamma'$. Indeed, because of the conventions for the sign of the
mean curvature, positive mean curvature means that the orientation
induced by the immersion $X$ is the same as the orientation induced by
the mean curvature vector $\vec{H}$ ; consequently, the basis
$(\eta,\Gamma',\vec{H})$ is indirect. We note
$\nu$ and $\eta$ the normal to $K$ and the conormal to $\Gamma$,
chosen as explained in section \ref{prel}.

The conormal $\eta$ to $\Gamma$ is a unit vector lying in the tangent
plane and normal to $\Gamma'(\tau)=\frac{\partial}{\partial\tau}X(\rho
e^{i\tau})$. Since the parametrization $X$ is conformal, the conormal
$\eta$ is necessarily colinear to $\frac{\partial}{\partial\rho}X(\rho
e^{i\tau})$. Since $\eta$ must point in the direction of $0\in\C$, we
have
$$\eta=-\frac{\frac{\partial}{\partial\rho}X(\rho e^{i\tau})}
{||\frac{\partial}{\partial\rho}X(\rho e^{i\tau})||}.$$

Then we have
\begin{eqnarray*}
\int_\Gamma<\eta,Y> & = & \int_0^{2\pi}<\eta,Y>
\left|\left|\frac{\partial}{\partial\tau}X\right|\right|\rmd\tau \\
& = & \int_0^{2\pi}<\eta,Y>
\left|\left|\frac{\partial}{\partial\rho}X\right|\right|\rho\rmd\tau \\
& = & \int_0^{2\pi}
-\rho\left<\frac{\partial}{\partial\rho}X,Y\right>\rmd\tau.
\end{eqnarray*}

According to lemma \ref{killingtrans}, we have
\begin{displaymath}
Y(\zeta,w)=\left(
\begin{array}{c}
-\frac{w^2}{\bar{\zeta}_0}
+\frac{(\zeta-\zeta_1)^2}{\zeta_0}-(\zeta-\zeta_1) \\
2w\re\frac{\zeta-\zeta_1}{\zeta_0}-w
\end{array}
\right).
\end{displaymath}

Consequently we have
$$\int_\Gamma<\eta,Y>=-\int_0^{2\pi}s_1(\rho,\tau)\rmd\tau$$ with
$$s_1(\rho,\tau)=
\frac{\rho}{w^2}\overline{\frac{\partial\zeta}{\partial\rho}}
\left(-\frac{w^2}{\bar{\zeta_0}}+\frac{(\zeta-\zeta_1)^2}{\zeta_0}
-(\zeta-\zeta_1)\right)
+\frac{\rho}{w}\frac{\partial w}{\partial\rho}
\left(2\frac{\zeta-\zeta_1}{\zeta_0}-1\right).$$

Let $\alpha$ be the canonical volume form of $\h^3$. We have 
$\alpha=\frac{1}{w^3}\rmd u\wedge\rmd v\wedge\rmd w$. Since $Y$ is a
Killing field,
we have $\mathrm{L}_Y\alpha=0$, so $0=\rmd(\mathrm{i}_Y\alpha)
+\mathrm{i}_Y\rmd\alpha=\rmd(\mathrm{i}_Y\alpha)$. Hence there exists
a 1-form $\beta$ such that $\mathrm{i}_Y\alpha=\rmd\beta$.

$\beta$ is the dual form of a vector field $Z$, {\it i.e.} we have
$\beta(\xi)=<Z,\xi>$. 

We compute that we can take
\begin{displaymath}
Z(\zeta,w)=\left(
\begin{array}{c}
i\frac{w^2}{\bar{\zeta}_0}\ln w+\frac{i}{2}\frac{(\zeta-\zeta_1)^2}{\zeta_0}
-\frac{i}{2}(\zeta-\zeta_1) \\
0
\end{array}
\right).
\end{displaymath}

Let $(e_1,e_2)$ be an orthonormal basis of the tangent space of $K$
such that the basis $(e_1,e_2,\nu)$ is direct. Then we have 
${\rm i}_Y\alpha(e_1,e_2)=\alpha(e_1,e_2,Y)=<\nu,Y>$. Consequently, on
$K$ the form ${\rm i}_Y\alpha$ is equal to $<\nu,Y>$ times the
canonical volume form of $K$. Hence we have 
$$\int_K<\nu,Y>=\int_K\mathrm{i}_Y\alpha.$$

On the other hand, Stokes's formula implies that 
$$\int_K\mathrm{i}_Y\alpha=-\int_\Gamma\beta$$
since we must take on $\Gamma$ the orientation given by $-\Gamma'$, as
explained before. 

Consequently we have $$\int_K<\nu,Y>=-\int_\Gamma\beta
=-\int_0^{2\pi}\left<\frac{\partial}{\partial\tau}X,Z\right>\rmd\tau
=-\int_0^{2\pi}\re(s_2(\rho,\tau))\rmd\tau$$
with
$$s_2(\rho,\tau)=\frac{1}{w^2}\overline{\frac{\partial\zeta}{\partial\tau}}
\left(i\frac{w^2}{\bar{\zeta}_0}\ln w
+\frac{i}{2}\frac{(\zeta-\zeta_1)^2}{\zeta_0}
-\frac{i}{2}(\zeta-\zeta_1)\right).$$

So $$\varphi=\int_\Gamma<\eta,Y>-2\int_K<\nu,Y>
=\int_0^{2\pi}\re(-s_1(\rho,\tau)+2s_2(\rho,\tau))\rmd\tau.$$

Since the real part does not change if we replace the first terms of
$s_1$ and $s_2$ by their conjugates, we obtain the expected result.
\end{proof}

\begin{lemma} \label{derivzeta}
%Let $E$ be a Bryant surface end, 
%$F=\left(
%\begin{array}{cc}
%A & B \\
%C & D
%\end{array} \right):\Omega^*\to\mathrm{SL}_2(\C)$
%its Bryant representation, $(\zeta,w):\Omega^*\to\h^3$ the
%corresponding immersion in the upper half-space model. Then we have
%$$\frac{1}{w^2}\overline{\left(\rho\frac{\partial\zeta}{\partial\rho}
%+i\frac{\partial\zeta}{\partial\tau}\right)}
%=2z(AB'-A'B).$$
We have the following identities :
\begin{equation} \label{relationab}
\frac{1}{w^2}\overline{\frac{\partial\zeta}{\partial\bar{z}}}
=AB'-A'B,
\end{equation}
\begin{equation} \label{relationw}
\frac{\rho}{w}\frac{\partial w}{\partial\rho}=
-z\frac{A'\bar{A}+B'\bar{B}}{|A|^2+|B|^2}
-\bar{z}\frac{A\bar{A}'+B\bar{B}'}{|A|^2+|B|^2},
\end{equation}
\begin{equation} \label{relationlog}
\frac{\partial}{\partial\tau}\ln{w}=
-iz\frac{A'\bar{A}+B'\bar{B}}{|A|^2+|B|^2}
+i\bar{z}\frac{A\bar{A}'+B\bar{B}'}{|A|^2+|B|^2}.
\end{equation}
\end{lemma}

\begin{proof}
Recall that $$\zeta=\frac{\bar{A}C+\bar{B}D}{|A|^2+|B|^2}$$ where $A$,
$B$, $C$ and $D$ are multivaluated holomorphic functions.
% Hence we
%have $\frac{\partial}{\partial\rho}A=e^{i\tau}A'$,
%$\frac{\partial}{\partial\rho}\bar{A}=e^{-i\tau}\bar{A}'$,
%$\frac{\partial}{\partial\tau}A=i\rho e^{i\tau}A'$,
%$\frac{\partial}{\partial\tau}\bar{A}=-i\rho e^{-i\tau}\bar{A}'$, and
%the analogous identities for $B$, $C$ and $D$.

We compute that
%$$\rho\frac{\partial\zeta}{\partial\rho}
%+i\frac{\partial\zeta}{\partial\tau}
%=2\rho e^{-i\tau}\frac{\bar{A}\bar{B}'-\bar{A}'\bar{B}}
%{(|A|^2+|B|^2)^2}.$$
\begin{eqnarray*}
\frac{\partial\zeta}{\partial\bar{z}} & = &
\frac{(\bar{A}'C+\bar{B}'D)(|A|^2+|B|^2)
-(\bar{A}C+\bar{B}D)(A\bar{A}'+B\bar{B}')}
{(|A|^2+|B|^2)^2} \\
& = & \frac{\bar{A}\bar{B}'-\bar{A}'\bar{B}}{(|A|^2+|B|^2)^2}
\end{eqnarray*}
because $AD-BC=1$.

And since $$w=\frac{1}{|A|^2+|B|^2},$$ we obtain relation
(\ref{relationab}).

Relations (\ref{relationw}) and (\ref{relationlog}) are consequences of elementary
computations using the fact that we have $\frac{\partial}{\partial\rho}A=e^{i\tau}A'$,
$\frac{\partial}{\partial\rho}\bar{A}=e^{-i\tau}\bar{A}'$,
$\frac{\partial}{\partial\tau}A=i\rho e^{i\tau}A'$,
$\frac{\partial}{\partial\tau}\bar{A}=-i\rho e^{-i\tau}\bar{A}'$, and
the analogous identities for $B$, $C$ and $D$ (because these are
multivaluated holomorphic functions).
\end{proof}

\begin{lemma} \label{coeffintegral}
We have $$s(\rho,\tau)=a_1(z)+\zeta_1a_2(z)+\frac{1}{\zeta_0}a_3(z)
+2\frac{\zeta_1}{\zeta_0}a_1(z)+\frac{\zeta_1^2}{\zeta_0}a_2(z)$$
where
$$a_1(z)=2z(B'C-A'D)+i\frac{\partial}{\partial\tau}\ln{w},$$
$$a_2(z)=2z(A'B-AB'),$$
$$a_3(z)=2z(C'D-CD')-2i\frac{\partial}{\partial\tau}(\zeta\ln{w})
+i\frac{\partial\zeta}{\partial\tau}.$$
\end{lemma}

\begin{proof}
We have the above expression for $s(\rho,\tau)$ with
$$a_1(z)=\frac{\zeta}{w^2}
\overline{\left(\rho\frac{\partial\zeta}{\partial\rho}
+i\frac{\partial\zeta}{\partial\tau}\right)}
+\frac{\rho}{w}\frac{\partial w}{\partial\rho},$$
$$a_2(z)=-\frac{1}{w^2}
\overline{\left(\rho\frac{\partial\zeta}{\partial\rho}
+i\frac{\partial\zeta}{\partial\tau}\right)},$$
$$a_3(z)=-\frac{\zeta^2}{w^2}
\overline{\left(\rho\frac{\partial\zeta}{\partial\rho}
+i\frac{\partial\zeta}{\partial\tau}\right)}
+\rho\frac{\partial\zeta}{\partial\rho}
-2i\frac{\partial\zeta}{\partial\tau}\ln{w}
-2\frac{\rho}{w}\frac{\partial w}{\partial\rho}\zeta.$$

The announced expression of $a_2(z)$ is a consequence of
formula (\ref{relationab}).

Because of formulae (\ref{relationab}) and (\ref{relationw}) we have
$$a_1(z)=2z(AB'-A'B)\zeta
-z\frac{A'\bar{A}+B'\bar{B}}{|A|^2+|B|^2}
-\bar{z}\frac{A\bar{A}'+B\bar{B}'}{|A|^2+|B|^2}.$$
Then a computation shows that
$$a_1(z)=2z(B'C-A'D)
+z\frac{A'\bar{A}+B'\bar{B}}{|A|^2+|B|^2}
-\bar{z}\frac{A\bar{A}'+B\bar{B}'}{|A|^2+|B|^2}.$$
Thus we obtain the above expression for $a_1(z)$ using formula (\ref{relationlog}).

Finally we have
\begin{eqnarray*}
a_3(z) & = & -a_1(z)\zeta
+\rho\frac{\partial\zeta}{\partial\rho}
-2i\frac{\partial}{\partial\tau}(\zeta\ln{w})
+2i\zeta\frac{\partial}{\partial\tau}\ln{w}
-\frac{\rho}{w}\frac{\partial w}{\partial\rho}\zeta \\
& = & -2z(B'C-A'D)\zeta
+\rho\frac{\partial\zeta}{\partial\rho}
-2i\frac{\partial}{\partial\tau}(\zeta\ln{w})
+i\zeta\frac{\partial}{\partial\tau}\ln{w}
-\frac{\rho}{w}\frac{\partial w}{\partial\rho}\zeta \\
& = & 2z(BC'-AD')\zeta
+2z\frac{\partial\zeta}{\partial z}
+i\frac{\partial\zeta}{\partial\tau}
-2i\frac{\partial}{\partial\tau}(\zeta\ln{w})
-2z\frac{\zeta}{w}\frac{\partial w}{\partial z} \\
& = & 2z(BC'-AD')\zeta
+2zw(\bar{A}C'+\bar{B}D')
+i\frac{\partial\zeta}{\partial\tau}
-2i\frac{\partial}{\partial\tau}(\zeta\ln{w}) \\
& = & 2z(C'D-CD')
+i\frac{\partial\zeta}{\partial\tau}
-2i\frac{\partial}{\partial\tau}(\zeta\ln{w}).
\end{eqnarray*}
\end{proof}

As an immediate consequence of lemmas \ref{integral} and \ref{coeffintegral} we obtain the following result.

\begin{lemma} \label{fluxzetazero}
%Let $E$ be a Bryant surface end given 
%by a conformal immersion $X=(\zeta,w):\Omega^*\to\h^3$ in
%the upper half-space model. Let
%$\zeta_0$ and $\zeta_1$ be two complex numbers, with $\zeta_0\neq 0$, 
%and let $Y$ be the Killing field associated to the translation along
%the geodesic $(\zeta_0+\zeta_1,\zeta_1)$.
Let $Y$ be the Killing field associated to the translation along
the geodesic $(\zeta_0+\zeta_1,\zeta_1)$.
Then the flux of $Y$ through $E$ is 
$$\varphi=\re\left(\varphi_1+\varphi_2\zeta_1
+\varphi_0\frac{1}{\zeta_0}
+2\varphi_1\frac{\zeta_1}{\zeta_0}
+\varphi_2\frac{\zeta_1^2}{\zeta_0}\right)$$
where $\varphi_0=4\pi\Res(D\rmd C-C\rmd D)$, $\varphi_1=4\pi\Res(C\rmd B-D\rmd A)$ and $\varphi_2=4\pi\Res(B\rmd A-A\rmd B)$.
\end{lemma}

Now we deal with the case where one of the extremities of the geodesic
is the point $\infty$.

\begin{lemma} \label{fluxinfinity}
Let $Y$ be the Killing field associated to the translation along the
geodesic $(\zeta_1,\infty)$. Then the flux of $Y$ through $E$ is
$$\varphi=\re(-\varphi_1-\varphi_2\zeta_1)$$
where $\varphi_0=4\pi\Res(D\rmd C-C\rmd D)$, $\varphi_1=4\pi\Res(C\rmd B-D\rmd A)$ and $\varphi_2=4\pi\Res(B\rmd A-A\rmd B)$.
\ref{fluxzetazero}.
\end{lemma}

\begin{proof}
We proceed as in lemmas \ref{integral}, \ref{coeffintegral} and
\ref{fluxzetazero}, replacing the expressions of $Y$ and $Z$ in lemma
\ref{integral} by
\begin{displaymath}
Y(\zeta,w)=\left(
\begin{array}{c}
\zeta-\zeta_1 \\
w
\end{array}
\right)
\end{displaymath}
and
\begin{displaymath}
Z(\zeta,w)=\left(
\begin{array}{c}
\frac{i}{2}(\zeta-\zeta_1) \\
0
\end{array}
\right).
\end{displaymath}
\end{proof}

\begin{thm} \label{fluxtranslation}
Let $\cC$ and $\cD$ be two distinct points in $\bar{\C}$. Let $Y$
be the Killing field associated to the translation along the 
geodesic $(\cC,\cD)$. Then the flux of $Y$ through $E$ is
$$\varphi=\re\left(\frac{\varphi_2\cC\cD+\varphi_1(\cC+\cD)+\varphi_0}
{\cC-\cD}\right)$$
where $\varphi_0=4\pi\Res(D\rmd C-C\rmd D)$, $\varphi_1=4\pi\Res(C\rmd B-D\rmd A)$ and $\varphi_2=4\pi\Res(B\rmd A-A\rmd B)$.
\end{thm}

\begin{proof}
If both $\cC$ and $\cD$ are different from $\infty$, then we set $\zeta_1=\cD$ and
$\zeta_0=\cC-\cD$, and the result comes from lemma \ref{fluxzetazero}.

If $\cD=\infty$ and $\cC\neq\infty$, then we set $\zeta_1=\cC$, and
the result comes from lemma \ref{fluxinfinity}.

If $\cC=\infty$ and $\cD\neq\infty$, then the result follows from the
above case and the fact that both the flux and the announced
expression are antisymmetric with respect to $(\cC,\cD)$.
\end{proof}

\subsection{Flux of Killing fields associated to rotations} \label{computationrotation}

\begin{lemma} \label{rotintegral}
Let $E$ be a Bryant surface end given 
by a conformal immersion $X=(\zeta,w):\Omega^*\to\h^3$ in
the upper half-space model. Let
$\zeta_0$ and $\zeta_1$ be two complex numbers, with $\zeta_0\neq 0$, 
and let $Y$ be the Killing field associated to the rotation about
the geodesic $(\zeta_0+\zeta_1,\zeta_1)$.
Then if $\rho$ is a sufficiently small positive number, the flux of
$Y$ through $E$ is 
$$\varphi=\int_0^{2\pi}\re(is(\rho,\tau))\rmd\tau$$
where $s(\rho,\tau)$ has been defined in lemma \ref{integral}.
\end{lemma}

\begin{proof}
We proceed as in lemma \ref{integral}, with 
\begin{displaymath}
Y(\zeta,w)=\left(
\begin{array}{c}
i\frac{w^2}{\bar{\zeta}_0}
+i\frac{(\zeta-\zeta_1)^2}{\zeta_0}-i(\zeta-\zeta_1) \\
-2w\im\frac{\zeta-\zeta_1}{\zeta_0}
\end{array}
\right)
\end{displaymath}
(see lemma \ref{killingrot}) and
\begin{displaymath}
Z(\zeta,w)=\left(
\begin{array}{c}
\frac{w^2}{\bar{\zeta}_0}\ln w-\frac{1}{2}\frac{(\zeta-\zeta_1)^2}{\zeta_0}
+\frac{1}{2}(\zeta-\zeta_1) \\
0
\end{array}
\right).
\end{displaymath}
\end{proof}

Using this lemma, we proceed as in section
\ref{computationtranslation} to compute the flux of Killing fields
associated to rotations.

\begin{thm} \label{fluxrotation}
Let $\cC$ and $\cD$ be two distinct points in $\bar{\C}$. Let $Y$
be the Killing field associated to the rotation about the 
geodesic $(\cC,\cD)$. Then the flux of $Y$ through $E$ is
$$\varphi=-\im\left(\frac{\varphi_2\cC\cD+\varphi_1(\cC+\cD)+\varphi_0}
{\cC-\cD}\right)$$
where $\varphi_0=4\pi\Res(D\rmd C-C\rmd D)$, $\varphi_1=4\pi\Res(C\rmd B-D\rmd A)$ and $\varphi_2=4\pi\Res(B\rmd A-A\rmd B)$.
\end{thm}

\subsection{Flux polynomial and equivalence with the residue-type flux matrix}

\begin{thm} \label{fluxpolynomial}
Let $E$ be a Bryant surface end 
whose Bryant representation is $F=\left(
\begin{array}{cc}
A & B \\
C & D
\end{array} \right):\Omega^*\to\mathrm{SL}_2(\C)$.
Then there exists a unique polynomial $P_E(X,Y)\in\C[X,Y]$ such
that, for all couples $(\cC,\cD)$ of distinct points in $\bar{\C}$, the flux
of the Killing field associated to the translation along the geodesic
$(\cC,\cD)$ through $E$ is $$\re\left(\frac{P_E(\cC,\cD)}{\cC-\cD}\right)$$
and the flux of the Killing field associated to the rotation about
the geodesic $(\cC,\cD)$ through $E$ is
$$-\im\left(\frac{P_E(\cC,\cD)}{\cC-\cD}\right)$$

This polynomial $P_E$ is symmetric and we have
$$P_E(X,Y)=\varphi_2XY+\varphi_1(X+Y)+\varphi_0$$ where
$\varphi_0=4\pi\Res(D\rmd C-C\rmd D)$,
$\varphi_1=4\pi\Res(C\rmd B-D\rmd A)$ and
$\varphi_2=4\pi\Res(B\rmd A-A\rmd B)$.

The polymomial $$\Pi_E(X)=P_E(X,X)=\varphi_2X^2+2\varphi_1X+\varphi_0$$
 is called the flux polynomial of $E$.
\end{thm}

\begin{proof}
This is a reformulation of theorems \ref{fluxtranslation} and
\ref{fluxrotation}.
\end{proof}

\begin{rem}
We have
$$\Pi_E(X)=-4\pi\Res\left(\omega^\#\left(X-\frac{1}{G}\right)^2\right).$$
\end{rem}

\begin{rem}
Knowing the flux polynomial is equivalent to knowing the flux of
Killing fields associated to all translations and rotations.
\end{rem}

In \cite{umeharaflux}, Rossman, Umehara and Yamada defined
a residue-type flux for Bryant surface ends. If an end $E$ is conformally
parametrized by $\Omega^*$ and has a Bryant representation 
$F=\left(
\begin{array}{cc}
A & B \\
C & D
\end{array} \right)$, then the flux matrix of $E$ is defined by
$$\Phi=-\frac{1}{2i\pi}\int_{\Gamma}(\rmd F)F^{-1}$$ where $\Gamma$
is a loop around $0$ with positive orientation. This matrix does not
depend on the choice of $\Gamma$. It is the residue at zero of the form 
$$-(\rmd F)F^{-1}=\left(
\begin{array}{cc}
C\rmd B-D\rmd A & B\rmd A-A\rmd B \\
C\rmd D-D\rmd C & B\rmd C-A\rmd D
\end{array} \right),$$
which is single-valued. Hence it does not depend on the parametrization.

Consequently, since $B\rmd C-A\rmd D=-(C\rmd B-D\rmd A)$, we have
$$\Phi=\frac{1}{4\pi}\left(
\begin{array}{cc}
\varphi_1 & \varphi_2 \\
-\varphi_0 & -\varphi_1
\end{array} \right).$$
Thus the coefficients of the flux matrix $\Phi$ are, up to constants, the same as the coefficients of the flux polynomial.

This proves the conjecture of Rossman, Umehara and Yamada (\cite{umeharaflux}, remark following example 8) : knowing the flux matrix $\Phi$ of the end $E$ is equivalent to knowing the flux through $E$ of all Killing fields associated to translations and rotations. 

We considered these two notions of flux for loops $\Gamma$ generating the fundamental group of an end. We can actually define these fluxes for any loop $\Gamma$ on a Bryant surface. We consider a neighbourhood of $\Gamma$ in the surface that is conformally parametrized by $\{z\in\C|1-\varepsilon<|z|<1+\varepsilon\}$ and such that $\Gamma$ is homologous to the curve corresponding to the circle $\{|z|=1\}$. Then the flux of a Killing field $Y$ through $\Gamma$ is equal to its flux through the curve corresponding to the circle $\{|z|=1\}$ (since the flux is a homology invariant). Thus we obtain, theorems \ref{fluxtranslation}, \ref{fluxrotation} and \ref{fluxpolynomial} with
$\varphi_0=-2i\int_{\{|z|=1\}}(D\rmd C-C\rmd D)$,
$\varphi_1=-2i\int_{\{|z|=1\}}(C\rmd B-D\rmd A)$ and
$\varphi_2=-2i\int_{\{|z|=1\}}(B\rmd A-A\rmd B)$.
These coefficients are, up to constants, the coefficients of the flux matrix $\Phi=-\frac{1}{2i\pi}\int_{\Gamma}(\rmd F)F^{-1}=
-\frac{1}{2i\pi}\int_{\{|z|=1\}}(\rmd F)F^{-1}$. Hence the two notions of flux are equivalent for any loop $\Gamma$ on the surface, and consequently for any homology class on the surface.

\begin{rem} It is easy to compute the flux matrix $\Phi$ of an end $E$ that is the image by a direct isometry of $\h^3$ of an end $E_0$ whose flux matrix $\Phi_0$ is known. Indeed, if $F$ and $F_0$ are the Bryant representations of $E$ and $E_0$, then there exists a matrix $P\in\mathrm{SL}_2(\C)$ such that $F_0=PF$. Then $\Phi=P^{-1}\Phi_0P$.
\end{rem}

\section{Embedded Bryant surface ends of finite total curvature}

Let us first recall and complete the results of S\'a Earp and Toubiana
(\cite{toubiana}).

Let $E$ be an enbedded Bryant surface end of finite total curvature
which is not part of a horosphere. We recall that $E$ is necessarily
regular (see \cite{yu}). Then, according to \cite{bryant},
the associated Weiertrass data have the following form :
\begin{eqnarray*} 
\left\{
\begin{array}{ccc}
g(z) & = & z^\mu f(z) \\
\omega & = & z^\nu h(z)\rmd z
\end{array} \right.
\end{eqnarray*}
in $\Omega^*$, where $f$ and $h$ are holomorphic
functions in a neighbourhood of zero such that $f(0)\neq 0$ and
$h(0)\neq 0$, and $\mu$ and $\nu$ are real numbers such that $\mu>0$,
$\nu\leqslant -1$, $\mu+\nu\in\Z$ and $\mu+\nu\geqslant-1$.

Since $f(0)\neq 0$, we can define a function $z\mapsto
f(z)^{\frac{1}{\mu}}$ in a neighbourhood of zero. Consequently,
we can replace $z$ by $zf(z)^{\frac{1}{\mu}}$ and assume that the
Weierstrass data have the following form :
\begin{eqnarray} \label{data}
\left\{
\begin{array}{ccc}
g(z) & = & z^\mu \\
\omega & = & z^\nu h(z)\rmd z.
\end{array} \right.
\end{eqnarray}

We distinguish two cases : the case where $\mu+\nu=-1$ will be dealt
with in section \ref{catends} and the case where $\mu+\nu\geqslant
0$ will be dealt with in section \ref{horoends}.

\subsection{Catenoidal ends} \label{catends}

\subsubsection{General representation} \label{genrep}

In this section we assume that $\mu+\nu=-1$. In this case the Hopf
differential $\omega\rmd g$ is of degree $-2$. Then, according to
\cite{toubiana}, we have $\mu\neq 1$ and, after replacing $f(z)$ by $1$,
\begin{equation} \label{h}
h(0)=\frac{1-\mu^2}{4\mu}
\end{equation}
and $$\frac{4\mu}{1-\mu}h'(0)=2\mu h'(0).$$ This second
equation implies that
\begin{equation} \label{hprime}
h'(0)=0.
\end{equation}

The Bryant representation of $E$ is given by 
\begin{eqnarray} \label{repbryant}
F=\left(
\begin{array}{cc}
A & B \\
C & D
\end{array} \right)
=\left(
\begin{array}{cc}
a_1z^{\lambda_1}f_1+a_2z^{\lambda_2}f_2 & 
b_1z^{r_1}r+b_2z^{r_2}g_2 \\
c_1z^{\lambda_1}f_1+c_2z^{\lambda_2}f_2 & 
d_1z^{r_1}r+d_2z^{r_2}g_2
\end{array} \right),
\end{eqnarray}
where $f_1$, $f_2$, $r$ and $g_2$ are holomorphic functions near 0 satifying
$f_1(0)=f_2(0)=r(0)=g_2(0)=1$, $\lambda_1=\frac{-1-\mu}{2}$,
$\lambda_2=\frac{1-\mu}{2}$,  $r_1=\frac{\mu-1}{2}$ and
$r_2=\frac{1+\mu}{2}$, and where $a_1$, $a_2$, $b_1$, $b_2$, $c_1$,
$c_2$, $d_1$ and $d_2$ are complex numbers satisfying
$a_1d_1-b_1c_1=0$, $a_1c_2-a_2c_1\neq 0$ and $b_1d_2-b_2d_1\neq 0$.

The functions $f_1$ and $f_2$ are such that 
$(z\mapsto z^{\lambda_1}f_1(z),z\mapsto z^{\lambda_2}f_2(z))$ is a
basis of the vector space of the solutions of the equation
$$X''-\frac{(z^{-1-\mu}h)'}{z^{-1-\mu}h}X'-\mu h z^{-2}X=0.$$ 
The functions $r$ and $g_2$ are such that 
$(z\mapsto z^{r_1}r(z),z\mapsto z^{r_2}g_2(z))$ is a
basis of the vector space of the solutions of the equation
$$X''-\frac{(z^{-1+\mu}h)'}{z^{-1+\mu}h}X'-\mu h z^{-2}X=0.$$ 

\begin{rem} \label{derivzero}
Since $\lambda_2=\lambda_1+1$, the function $f_2$ is uniquely defined, and
the function $f_1$ is uniquely defined if we fix the value of its
derivative at zero. In the same way, since
$r_2=r_1+1$, the function $g_2$ is uniquely defined, and
the function $r$ is uniquely defined if we fix the value of its
derivative at zero.
\end{rem}

From the identity $\omega=A\rmd C-C\rmd A$ (see \cite{umehara} or
\cite{rosenberg}) 
we obtain that
\begin{equation} \label{relomega}
h=(a_1c_2-a_2c_1)(f_1f_2-zf_1'f_2+zf_1f_2').
\end{equation}

Taking the order 1 terms, we get
\begin{equation}
f_2'(0)=0.
\end{equation}

In the same way, from the identity $g^2\omega=B\rmd D-D\rmd B$ (see
\cite{umehara} or \cite{rosenberg}) 
we obtain that
\begin{equation} \label{relgomega}
h=(b_1d_2-b_2d_1)(rg_2-zr'g_2+zrg_2').
\end{equation}

Taking the order 1 terms, we get
\begin{equation}
g_2'(0)=0.
\end{equation}

\subsubsection{Canonical representation} \label{canrep}

Sá Earp and Toubiana (\cite{toubiana}) have shown that we can reduce
ourselves to a more simple Bryant representation up to an isometry of
$\h^3$. More precisely, we can choose complex
numbers $\alpha$, $\beta$, $\gamma$ and $\delta$ satisfying
$\alpha\delta-\beta\gamma=1$, $\alpha a_1+\beta c_1=\alpha b_1+\beta d_1
=\gamma a_2+\delta c_2=0$, and $\alpha a_2+\beta c_2=1$.
If we replace $F=\left(
\begin{array}{cc}
A & B \\
C & D
\end{array} \right)$ by
$\left(
\begin{array}{cc}
\alpha & \beta \\
\gamma & \delta
\end{array} \right)
\left(
\begin{array}{cc}
A & B \\
C & D
\end{array} \right)$, we obtain an end which is the image of $E$ by a direct
isometry $\Psi$ of $\h^3$, which has the same Weierstrass data as $E$,
and whose Bryant representation is given by
\begin{eqnarray} \label{repbryantcanon}
\left(
\begin{array}{cc}
A(z) & B(z) \\
C(z) & D(z)
\end{array} \right)=\left(
\begin{array}{cc}
z^{\lambda_2}f_2(z) & \frac{\mu-1}{\mu+1}z^{r_2}g_2(z) \\
\frac{\mu^2-1}{4\mu}z^{\lambda_1}f_1(z) &
\frac{(1+\mu)^2}{4\mu}z^{r_1}g_1(z)
\end{array} \right),
\end{eqnarray}
where $g_1$ is a holomorphic function near 0 satifying
$g_1(0)=1$. The isometry $\Psi$ induces on
$\bar{\C}$ the map 
$\zeta\mapsto\frac{\delta\zeta+\gamma}{\beta\zeta+\alpha}$, which we
also denote $\Psi$.

\begin{defn} \label{canonical}
Let $\mu\in(0,1)\cup(1,\infty)$ and ${\mathcal Z}\in\C$. An end which has 
Weierstrass data given by (\ref{data}) and a Bryant
representation $F$ given by (\ref{repbryantcanon}), where $f_1$ has
been chosen such that $${\mathcal Z}=\frac{\mu^2-1}{4\mu}f_1'(0)$$
(see remark \ref{derivzero}),
is called a canonical catenoidal end of growth $1-\mu$, 
of asymptotic boundary $\infty$ and of axis $({\mathcal Z},\infty)$, 
and the Weierstrass data given by (\ref{data}) and
the Bryant representation $F$ given by (\ref{repbryantcanon}) are 
called respectively its canonical Weierstrass representation and its
canonical Bryant representation. 
\end{defn}

We now explain this terminology by giving a geometric description of
such an end.

\begin{prop} \label{description}
Let $\mu\in(0,1)\cup(1,\infty)$, ${\mathcal Z}\in\C$ and $E$ be a
canonical cate\-noidal end of growth $1-\mu$,  
of asymptotic boundary $\infty$ and of axis $({\mathcal Z},\infty)$. Then
there exists a parametrization $(w,\tau)\mapsto(\zeta(w,\tau),w)$ of
$E$ in the upper half-space model of $\h^3$, and a parametrization
$(w,\tau)\mapsto(\tilde{\zeta}(w,\tau),w)$ of a half-catenoid cousin
of growth $1-\mu$, 
of asymptotic boundary $\infty$ and of axis $(0,\infty)$,
such that
$$\zeta(w,\tau)=\tilde{\zeta}(w,\tau)+\cZ+\rmo(1)$$ when $w$ tends to
$\infty$ if $\mu<1$ and to $0$ if $\mu>1$.
\end{prop}

\begin{proof}
We assume that the Weierstrass data of $E$ are given by (\ref{data})
and its Bryant representation $F$ by (\ref{repbryantcanon}), with
${\mathcal Z}=\frac{\mu^2-1}{4\mu}f_1'(0)$ (see remark \ref{derivzero}).

Because of formulae (\ref{zetaabcd}) and (\ref{wabcd}), in the upper
half-space model the end $E$ is given by 
$$\zeta(z)=(u+iv)(z)=\frac{\mu^2-1}{4\mu z}
\frac{f_1\bar{f_2}+g_1\bar{g_2}|z|^{2\mu}}
{|f_2|^2+(\frac{\mu-1}{\mu+1})^2|g_2|^2|z|^{2\mu}},$$
$$w(z)=\frac{|z|^{\mu-1}}{|f_2|^2+(\frac{\mu-1}{\mu+1})^2|g_2|^2|z|^{2\mu}}.$$

From this we deduce that the asymptotic boundary of $E$ is actually $\infty$.

Define $$\tilde{\zeta}(z)=\frac{\mu^2-1}{4\mu z}
\frac{1+|z|^{2\mu}}{1+(\frac{\mu-1}{\mu+1})^2|z|^{2\mu}}$$ and
$$\tilde{w}(z)=\frac{|z|^{\mu-1}}{1+(\frac{\mu-1}{\mu+1})^2|z|^{2\mu}}.$$
$\tilde{\zeta}$ and $\tilde{w}$ are the coordinates of the catenoid
cousin of growth $1-\mu$ and of
axis of revolution $(0,\infty)$, such that the end at $\infty\in\bar{\C}$
corresponds to $z=0$. $\tilde{w}$ depends only on $|z|$.

We have $\zeta(z)=\tilde{\zeta}(z)+{\mathcal Z}+\rmo(1)$ since ${\mathcal
  Z}=\frac{\mu^2-1}{4\mu}f_1'(0)$, and
$w(z)=\tilde{w}(z)(1+\rmO(z^2))$ since $f_2'(0)=g_2'(0)=0$. 

Let $(\rho,\tau)$ denote the polar coordinates in $\Omega$ ({\it i.e.}
$z=\rho e^{i\tau}$).
Since $\frac{\partial w}{\partial\rho}\neq 0$ for $\rho$
sufficiently small, we can do the change of parameters
$(\rho,\tau)\mapsto(w,\tau)$. 

Since $f_2'(0)=g_2'(0)=0$, we have the following asymptotic development :
$$w=\rho^{\mu-1}\left(1+\sum_{j=1}^{p}\alpha_j\rho^{2j\mu}
+\rmO(\rho^2)\right)$$ where $p$ is the largest integer such that
$2p\mu<2$ and the $\alpha_j$ are real constants which depend only on
$\mu$.

Consequently, we have the following asymptotic development for the
inverse function :
$$\rho=w^{\frac{1}{\mu-1}}\left(1+
\sum_{j=1}^{p}\beta_jw^{\frac{2j\mu}{\mu-1}}
+\rmO(w^{\frac{2}{\mu-1}})\right)$$ when $w$ tends to $\infty$ if
$\mu<1$ and to $0$ if $\mu>1$, and where the $\beta_j$ are real
constants which depend only on $\mu$.

We also have $$\zeta=\frac{\mu^2-1}{4\mu\rho e^{i\tau}}
\left(1+\sum_{j=1}^{q}\gamma_j\rho^{2j\mu}\right)+\cZ+\rmo(1)$$ where
$q$ is the largest integer such that $2q\mu\leqslant 1$ and the $\gamma_j$ are
real constants which depend only on $\mu$.

Reporting the asymptotic development of $w$, we get
\begin{eqnarray*}
\zeta(w,\tau) & = & \frac{\mu^2-1}{4\mu e^{i\tau}}w^{-\frac{1}{\mu-1}}
\left(1+\sum_{j=1}^{p}\delta_jw^{\frac{2j\mu}{\mu-1}}
+\rmO(w^{\frac{2}{\mu-1}})\right)+\cZ+\rmo(1) \\
& = & \frac{\mu^2-1}{4\mu e^{i\tau}}w^{-\frac{1}{\mu-1}}
\left(1+\sum_{j=1}^{q}\delta_jw^{\frac{2j\mu}{\mu-1}}
\right)+\cZ+\rmo(1)
\end{eqnarray*}
where the $\delta_j$ are real constants which depend only on $\mu$.

The same arguments hold for the canonical catenoid of axis $(0,\infty)$
parame\-trized by $(\tilde{\zeta},\tilde{w})$. Consequently we get
$$\zeta(w,\tau)=\tilde{\zeta}(w,\tau)+\cZ+\rmo(1).$$ 
\end{proof}

This means that
the end $E$ is asymptotic, in the neighbourhood of $\infty\in\bar{\C}$,
to a half-catenoid cousin of growth $1-\mu$ and of
axis of revolution $({\mathcal Z},\infty)$, in a stronger sense than the
sense defined in \cite{toubiana} (in \cite{toubiana}, two half-catenoid cousins
whose asymptotic boundary is $\infty$ and having the same growth
are asymptotic to each other up to a Euclidean homothety,
independently of their axes). The complex number $\cZ$ is the only one
with this property. 

\begin{defn} \label{catend}
Let $\mu\in(0,1)\cup(1,\infty)$ and ${\mathcal A}$, ${\mathcal B}$ two
distinct points in $\bar{\C}$. Let
$E$ be an embedded Bryant surface end of finite total curvature
which is not part of a horosphere.
We say that $E$ is a catenoidal end of growth $1-\mu$, 
of asymptotic boundary ${\mathcal B}$ and of axis $({\mathcal A},{\mathcal B})$ if
there exists an isometry of $\h^3$ (direct or not)
which maps ${\mathcal A}$ to $0$, ${\mathcal B}$ to $\infty$ and $E$ to a canonical
catenoidal end of growth $1-\mu$,
of aymptotic boundary $\infty$ and of axis $(0,\infty)$.
\end{defn}

A half-catenoid cousin of growth $1-\mu$, of aymptotic boundary $\cB$
and of axis of revolution $(\cA,\cB)$ is of course a catenoidal
end of growth $1-\mu$, of aymptotic boundary $\cB$ and of axis
$(\cA,\cB)$. 

A canonical catenoidal end of growth $1-\mu$, of aymptotic boundary
$\infty$ and of axis $(\cZ,\infty)$ is a catenoidal end of growth
$1-\mu$, of aymptotic boundary $\infty$ and of axis $(\cZ,\infty)$ :
it suffices to consider the isometry $(\zeta,w)\mapsto(\zeta-{\mathcal Z},w)$.

We can now prove the following theorem.

\begin{thm} \label{uniqueness}
Let $E$ be an embedded Bryant surface end of finite total curvature
which is not part of a horosphere. Assume that its Weierstrass data
are given by (\ref{data}) with $\mu+\nu=-1$. Then
there exist a unique real $\chi$ and a unique couple of distinct
points $(\cA,\cB)$  
such that $E$ is a catenoidal end of growth $\chi$, 
of asymptotic boundary ${\mathcal B}$ and of axis 
$({\mathcal A},{\mathcal B})$. 

Moreover we have $\chi=1-\mu$ and, if the Bryant representation of
$E$ is given by 
\begin{eqnarray*}
F(z)=\left(
\begin{array}{cc}
z^{\lambda_1}a(z) & z^{r_1}b(z) \\
z^{\lambda_1}c(z) & z^{r_1}d(z)
\end{array} \right),
\end{eqnarray*}
we have $\cA=c'(0)/a'(0)$ and $\cB=c(0)/a(0)$.
\end{thm}

\begin{proof} The existence has already been proved in section
\ref{genrep} and in the beginning of section \ref{canrep} (choosing
$f_1'(0)=0$, see remark \ref{derivzero}). 

The uniqueness of $\cB$ is clear, since the asymptotic boundary of $E$
is the set of its accumulation points in $\bar{\C}$.

Assume that there exist two points $\cA_1$ and $\cA_2$ and two numbers
$\mu_1$ and $\mu_2$ such that $E$
is both a catenoidal end of growth $1-\mu_1$ and of axis $(\cA_1,\cB)$
and a catenoidal end of growth $1-\mu_2$ and of axis
$(\cA_2,\cB)$. Then there exists an isometry $\Psi_1$ of $\h^3$ which
maps $\cA_1$ to $0$, $\cB$ to $\infty$ and $E$ to a canonical
catenoidal end of growth $1-\mu_1$ and of axis $(0,\infty)$, and there
exists an isometry $\Psi_2$ of $\h^3$ which 
maps $\cA_2$ to $0$, $\cB$ to $\infty$ and $E$ to a canonical
catenoidal end of growth $1-\mu_2$ and of axis $(0,\infty)$.

Consequently there exists a parametrization
$(w,\tau)\mapsto(\zeta_2(w,\tau),w)$ of the end $\Psi_2(E)$ such that
$\zeta_2(w,\tau)=\tilde{\zeta}_{\mu_2}(w,\tau)+\rmo(1)$ when $w$ tends
to $\infty$ if $\mu_2<1$ and to $0$ if $\mu_2>1$, where
$\tilde{\zeta}_{\mu_2}$ corresponds to the canonical catenoid of
growth $1-\mu_2$ and of axis $(0,\infty)$.

The isometry $\Psi_1\circ\Psi_2^{-1}$ fixes $\infty$ and maps $0$ to
$\cZ=\Psi_1(\cA_2)$. Assume that this isometry is direct. Then it is
the composition of a twist about $(0,\infty)$ and of the Euclidean
translation by the vector $\cZ$. Consequently, the end $\Psi_1(E)$ has a
parametrisation of the form
$(w,\tau)\mapsto(\lambda\zeta_2(w/|\lambda|,\tau)+\cZ,w)$ with
$\lambda\in\C^*$. 

On the other hand, there exists a parametrization
$(w,\tau')\mapsto(\zeta_1(w,\tau'),w)$ of $\Psi_1(E)$ such that
$\zeta_1(w,\tau')=\tilde{\zeta}_{\mu_1}(w,\tau')+\rmo(1)$ when $w$ tends
to $\infty$ if $\mu_1<1$ and to $0$ if $\mu_1>1$, where
$\tilde{\zeta}_{\mu_1}$ corresponds to the canonical catenoid of
growth $1-\mu_1$ and of axis $(0,\infty)$.

The numbers $1-\mu_1$ and $1-\mu_2$ must have the same sign, since $w$
cannot tend to both $\infty$ and $0$. And there
exists $\tau$ such that $\cZ=c(w)\tilde{\zeta}_{\mu_1}(w,\tau)$
with $c(w)\geqslant 0$ and for each $w$ there exists $\tau'(w)$ such that 
$\zeta_1(w,\tau)=\lambda\zeta_2(w/|\lambda|,\tau'(w))$. We deduce that
$$\lambda\tilde{\zeta}_{\mu_2}(w/|\lambda|,\tau'(w))+\cZ
-\tilde{\zeta}_{\mu_1}(w,\tau)=\rmo(1).$$ Using the expressions of
$\tilde{\zeta}_{\mu_1}$ and $\tilde{\zeta}_{\mu_2}$, we obtain that
$\mu_1=\mu_2$, $e^{i(\tau'(w)-\tau)}\to\lambda$ and hence
$|\lambda|=1$. Writing $\lambda=e^{i\theta}$, we have
$$(e^{i(\theta+\tau-\tau'(w))}+c(w)-1)
\tilde{\zeta}_{\mu_1}(w,\tau)=\rmo(1),$$ and so
\begin{equation} \label{axe}
e^{i(\theta+\tau-\tau'(w))}+c(w)-1=
\rmo(\tilde{\zeta}_{\mu_1}(w,\tau)^{-1}).
\end{equation}

Taking the imaginary part in (\ref{axe}) we get
$$\sin(\theta+\tau-\tau'(w))=\rmo(\tilde{\zeta}_{\mu_1}(w,\tau)^{-1}),$$
and consequently 
$$\cos(\theta+\tau-\tau'(w))-1=\rmo(\tilde{\zeta}_{\mu_1}(w,\tau)^{-1}).$$
On the other hand, taking the real part in (\ref{axe}) we get
$$\cos(\theta+\tau-\tau'(w))+c(w)-1
=\rmo(\tilde{\zeta}_{\mu_1}(w,\tau)^{-1}),$$
so $$c(w)=\rmo(\tilde{\zeta}_{\mu_1}(w,\tau)^{-1}),$$ and finally
$$\cZ=c(w)\tilde{\zeta}_{\mu_1}(w,\tau)=\rmo(1).$$ This means that
$\cZ=0$. We conclude that $\cA_2=\cA_1$.  

If the isometry $\Psi_1\circ\Psi_2^{-1}$ is indirect, then it is the
composition of the symmetry about the plane $\{\re\zeta=0\}$ and of the two
aforementioned isometries, so the same arguments hold, replacing
$\tilde{\zeta}_{\mu_2}$ by its conjugate.

To complete the proof, it now suffices to compute the values of $\cA$
and $\cB$. Using the notations of the beginning of section
\ref{canrep} with $f_1'(0)=0$ (see remark \ref{derivzero}), we have
$a=a_1f_1+a_2zf_2$ and $c=c_1f_1+c_2zf_2$. Hence we get $a(0)=a_1$,
$a'(0)=a_2$, $c(0)=c_1$ and $c'(0)=c_2$. The expression of $\cB$
follows from formulae (\ref{zetaabcd}) and (\ref{wabcd}). And since
$\gamma a_2+\delta c_2=0$, we have $\Psi(c_2/a_2)=0$ (even if
$a_2=0$), so $\cA=c_2/a_2=c'(0)/a'(0)$.
\end{proof}

\begin{rem}
The fact that $|\lambda|=1$ means that among all the half-catenoid
cousins of growth $1-\mu$, of asymptotic boundary $\cB$ and of axis
$(\cA,\cB)$ there exists a unique one to which $E$ is strongly
asymptotic.
\end{rem}

The following fact is now clear.

\begin{prop}
Let $E$ be a catenoidal end of growth $1-\mu$, 
of asymptotic boundary ${\mathcal B}$ and of axis $(\cA,\cB)$. Let $\Psi$
be an isometry of $\h^3$ (direct or not). Then $\Psi(E)$ is a
catenoidal end of growth $1-\mu$, of asymptotic boundary $\Psi(\cB)$
and of axis $(\Psi(\cA),\Psi(\cB))$.
\end{prop}

\begin{rem}
In definition \ref{catend} we can require the isometry to be direct.
\end{rem}

\begin{rem}
The notion of canonical end has no geometrical meaning, but it will be
more convenient to use this terminology to compute the flux (see
section \ref{sectionfluxcat}). Any catenoidal end of axis
$(\cZ,\infty)$ is the image of a canonical one by a twist about
$(\cZ,\infty)$.
\end{rem}

\subsection{Horospherical ends} \label{horoends}

\subsubsection{General representation}

In this section, we assume that $E$ is an end whose Weierstrass
data are given by (\ref{data}) with $\mu+\nu\geqslant 0$. Since we have
a single-valued embedding, according to \cite{toubiana} we have
$\nu=-2$, $\mu\in\N$, $\mu\geqslant 2$, 
\begin{eqnarray} \label{horohprime}
\left\{
\begin{array}{ccc}
h'(0)=2h(0)^2 & \mathrm{if} & \mu=2, \\
h'(0)=0 & \mathrm{if} & \mu\geqslant 3.
\end{array} \right.
\end{eqnarray}

The Bryant representation of $E$ is given by 
\begin{eqnarray} \label{hrepbryant}
F=\left(
\begin{array}{cc}
A & B \\
C & D
\end{array} \right)=\left(
\begin{array}{cc}
a_1z^{-1}f_1+a_2f_2 & 
b_1r+b_2z^{2\mu-1}g_2 \\
c_1z^{-1}f_1+c_2f_2 & 
d_1r+d_2z^{2\mu-1}g_2
\end{array} \right),
\end{eqnarray}
where $f_1$, $f_2$, $r$ and $g_2$ are holomorphic functions near 0 satifying
$f_1(0)=f_2(0)=r(0)=g_2(0)=1$, and where $a_1$, $a_2$, $b_1$, $b_2$, $c_1$,
$c_2$, $d_1$ and $d_2$ are complex numbers satisfying
$a_1d_1-b_1c_1=0$, $a_1c_2-a_2c_1\neq 0$ and $b_1d_2-b_2d_1\neq 0$.

The functions $f_1$ and $f_2$ are such that 
$(z\mapsto z^{-1}f_1(z),z\mapsto f_2(z))$ is a
basis of the vector space of the solutions of the equation
$$X''-\frac{(z^{-2}h)'}{(z^{-2}h)}X'-\mu h z^{\mu-3}X=0.$$ 
The functions $r$ and $g_2$ are such that 
$(z\mapsto r(z),z\mapsto z^{2\mu-1}g_2(z))$ are a
basis of the vector space of the solutions of the equation
$$X''-\frac{(z^{2\mu-2}h)'}{(z^{2\mu-2}h)}X'-\mu h z^{\mu-3}X=0.$$ 

\begin{rem} \label{derivzerohoro}
The function $f_2$ is uniquely defined, and
the function $f_1$ is uni\-quely defined if we fix the value of its
derivative at zero.
\end{rem}

\subsubsection{Canonical representation} \label{horocanrep}

As for catenoidal ends, Sá Earp and Toubia\-na (\cite{toubiana}) have
shown that we can reduce ourselves to a more simple Bryant
representation up to an isometry of $\h^3$. More precisely, we can
choose complex numbers $\alpha$, $\beta$, $\gamma$ and $\delta$ satisfying
$\alpha\delta-\beta\gamma=1$, $\alpha a_1+\beta c_1=\alpha b_1+\beta d_1
=\gamma a_2+\delta c_2=0$, and $\alpha a_2+\beta c_2=1$.
If we replace $F=\left(
\begin{array}{cc}
A & B \\
C & D
\end{array} \right)$ by
$\left(
\begin{array}{cc}
\alpha & \beta \\
\gamma & \delta
\end{array} \right)
\left(
\begin{array}{cc}
A & B \\
C & D
\end{array} \right)$, we obtain an end which is the image of $E$ by a direct
isometry $\Psi$ of $\h^3$, which has the same Weierstrass data as $E$,
and whose Bryant representation is given by
\begin{eqnarray} \label{hororepbryantcanon}
F(z)=\left(
\begin{array}{cc}
A(z) & B(z) \\
C(z) & D(z)
\end{array} \right)=\left(
\begin{array}{cc}
f_2(z) & bz^{2\mu-1}g_2(z) \\
cz^{-1}f_1(z) & g_1(z)
\end{array} \right),
\end{eqnarray}
where $g_1$ is a holomorphic function near 0 satifying $g_1(0)=1$,
$b\in\C^*$ and $c\in\C^*$. 
The isometry $\Psi$ induces on
$\bar{\C}$ the map 
$\zeta\mapsto\frac{\delta\zeta+\gamma}{\beta\zeta+\alpha}$, which we
also denote $\Psi$.

\begin{defn} \label{horocanonical}
An end which has 
Weierstrass data given by (\ref{data}) and a Bryant
representation $F$ given by (\ref{hororepbryantcanon})
is called a canonical horospherical end of asymptotic boundary $\infty$, 
and the Weierstrass data given by (\ref{data}) and
the Bryant representation $F$ given by (\ref{repbryantcanon}) are 
called respectively its canonical Weierstrass representation and its
canonical Bryant representation. 
\end{defn}

We now assume that the end $E$ has Weierstrass data given by (\ref{data}) and a
Bryant representation $F$ given by (\ref{hororepbryantcanon}). 

Because of formulae (\ref{zetaabcd}) and (\ref{wabcd}), in the upper
half-space model the end $E$ is given by 
$$\zeta(z)=(u+iv)(z)=\frac{c}{z}
\frac{f_1\bar{f_2}+\frac{b}{c}z\bar{z}^{2\mu-1}g_1\bar{g_2}}
{|f_2|^2+|b|^2|z|^{2\mu}|g_2|^2},$$
$$w(z)=\frac{1}{|f_2|^2+|b|^2|z|^{2\mu}|g_2|^2}.$$

From the identity $\omega=A\rmd C-C\rmd A$ we obtain that
$$c(-f_1f_2+zf_1'f_2-zf_1f_2')=h.$$
Taking the order zero term, we get
\begin{equation} \label{horoc}
c=-h(0).
\end{equation}
Taking the order one term, we get
\begin{equation} \label{horof2prime}
h'(0)=-2cf_2'(0).
\end{equation}

Taking the order one term in the identity $AD-BC=1$, we get 
\begin{equation} \label{horofg}
f_2'(0)+g_1'(0)=0.
\end{equation}

Since the end has finite total curvature and is regular, we can write
$$\omega\rmd g=\sum_{j=-2}^{\infty}q_jz^j\rmd z^2.$$

We compute that $$\omega\rmd g=\mu z^{\mu-3}h(z)\rmd z^2.$$
Hence we have $$q_{-2}=0$$ and
\begin{eqnarray} \label{horohopf}
\left\{
\begin{array}{ccc}
q_{-1}=2h(0) & \mathrm{if} & \mu=2, \\
q_{-1}=0 & \mathrm{if} & \mu\geqslant 3.
\end{array} \right.
\end{eqnarray}

\subsection{Classification}

Here we summarize the results we have obtained.

\begin{thm} \label{classification}
Let $E$ be an embedded Bryant surface end of finite total
curvature. Then we are in one of the following cases :
\begin{itemize}
{\item $E$ is part of a horosphere,}
{\item $E$ is not part of a horosphere and there exists a point 
${\mathcal B}\in\bar{\C}$ such that $E$ is a horospherical end of
asymptotic boundary ${\mathcal B}$,}
{\item $E$ is not part of a horosphere and there exist a real 
$\mu\in(0,1)\cup(1,\infty)$ and two distinct points 
${\mathcal A},{\mathcal B}\in\bar{\C}$ such that $E$ is a catenoidal end of
growth $1-\mu$, of asymptotic boundary ${\mathcal B}$, and of axis 
$({\mathcal A},{\mathcal B})$.}
\end{itemize}
\end{thm}

\begin{proof}
It sufficies to show that we cannot be in two cases at the same
time. This is a consequence of the fact that the Hopf differential is
zero for horospheres, is non-zero and has a degree greater than or
equal to $-1$ for horospherical ends, and has a degree equal to $-2$
for catenoidal ends.
\end{proof}

\section{Flux for embedded ends of finite total curvature}

\subsection{Flux for catenoidal ends} \label{sectionfluxcat}

\begin{lemma} \label{fluxcatends}
Let $\mu\in(0,1)\cup(1,\infty)$ and ${\mathcal Z}\in\C$. Let $E$ be a canonical
catenoidal end of growth 
$1-\mu$, of asymptotic boundary $\infty$ and of axis $({\mathcal Z},\infty)$. Let
$\zeta_0$ and $\zeta_1$ be two complex numbers, with $\zeta_0\neq 0$.
Then the flux polynomial of $E$ is
$$\Pi_E(X)=2\pi(\mu^2-1)(X-\cZ),$$
the flux of the Killing field associated to the translation along
the geodesic $(\zeta_0+\zeta_1,\zeta_1)$ through $E$ is 
$$\pi(\mu^2-1)\left(2\re\left(\frac{\zeta_1-{\mathcal Z}}{\zeta_0}\right)+1\right),$$
the flux of the Killing field associated to the rotation about
the geodesic $(\zeta_0+\zeta_1,\zeta_1)$ through $E$ is 
$$2\pi(1-\mu^2)\im\left(\frac{\zeta_1-{\mathcal Z}}{\zeta_0}\right),$$
the flux of the Killing field associated to the translation along
the geodesic $(\zeta_1,\infty)$ through $E$ is 
$$\pi(1-\mu^2),$$
and the flux of the Killing field associated to the rotation about
the geodesic $(\zeta_1,\infty)$ through $E$ is zero.
\end{lemma}

\begin{proof}
Using the canonical Bryant representation (\ref{repbryantcanon}), we compute that the coefficients of the flux polynomial are
$$\varphi_0=4\pi\Res(D\rmd C-C\rmd D)=2\pi(1-\mu^2)\cZ,$$
$$\varphi_1=4\pi\Res(C\rmd B-D\rmd A)=\pi(\mu^2-1),$$
$$\varphi_2=4\pi\Res(B\rmd A-A\rmd B)=0.$$
Applying theorems \ref{fluxtranslation}, \ref{fluxrotation} and \ref{fluxpolynomial}, we obtain the announced results.
\end{proof}

\begin{thm} \label{theoremfluxcatenoidal}
Let $\mu\in(0,1)\cup(1,\infty)$.
Let ${\mathcal A}$, ${\mathcal B}$, ${\mathcal C}$ and ${\mathcal D}$ be four points
in $\bar{\C}$ 
such that ${\mathcal A}\neq{\mathcal B}$ 
and ${\mathcal C}\neq{\mathcal D}$. 
%Let $Y$ be the Killing field associated to the
%translation along the geodesic $({\mathcal C},{\mathcal D})$.
Let $E$ be a
catenoidal end of growth $1-\mu$, of asymptotic boundary ${\mathcal B}$
and of axis $({\mathcal A},{\mathcal B})$.
Then the flux of the Killing field associated to the translation along the geodesic $({\mathcal C},{\mathcal D})$ through $E$ is
$$\pi(1-\mu^2)(2\re(\cA,\cC,\cD,\cB)-1),$$
Then the flux of the Killing field associated to the rotation about the geodesic $({\mathcal C},{\mathcal D})$ through $E$ is
$$-2\pi(1-\mu^2)\im(\cA,\cC,\cD,\cB),$$
and the flux polynomial of $E$ is
$$\Pi_E(X)=2\pi(1-\mu^2)\frac{(X-\cA)(X-\cB)}{\cB-\cA}.$$
(In the case where $\cA=\infty$ (respectively $\cB=\infty$), the above
formula means $\Pi_E(X)=2\pi(1-\mu^2)(X-\cB)$ (respectively
$\Pi_E(X)=-2\pi(1-\mu^2)(X-\cA)$).)
\end{thm}

\begin{proof}
We first compute the flux $\varphi$ of the Killing field $Y$ associated to the
translation along the geodesic $({\mathcal C},{\mathcal D})$.

According to what has been done in section \ref{catends}, the end $E$ has
Weierstrass data given by (\ref{data}), a Bryant representation $F$
given by (\ref{repbryant}), and, given a complex number ${\mathcal Z}$,
there exists a direct isometry $\Psi$ of $\h^3$ which maps $E$ to a canonical
catenoidal end of growth $1-\mu$, of asymptotic boundary $\infty$ and
of axis $({\mathcal Z},\infty)$.

Assume that neither $\cC$ nor $\cD$ is equal to $\cB$. Set
$\zeta_0=\Psi({\mathcal C})-\Psi({\mathcal D})$ and $\zeta_1=\Psi({\mathcal D})$.
Then $\zeta_0$ and $\zeta_1$ are different from $\infty$, and $\Psi$
maps ${\mathcal A}$ to ${\mathcal Z}$, ${\mathcal B}$ to $\infty$, 
${\mathcal C}$ to $\zeta_0+\zeta_1$ and ${\mathcal D}$ to $\zeta_1$. Hence the
flux $\varphi$ of $Y$ through $E$ is equal to the flux of the Killing
field associated to the 
translation along the geodesic $(\zeta_0+\zeta_1,\zeta_1)$ through
$\Psi(E)$. This flux has been calculated in lemma
\ref{fluxcatends} : we have
$$\varphi=
\pi(\mu^2-1)\left(2\re\left(\frac{\zeta_1-{\mathcal Z}}{\zeta_0}\right)+1\right).$$

We compute that $$({\mathcal Z},\zeta_0+\zeta_1,\infty,\zeta_1)
=-\frac{\zeta_1-{\mathcal Z}}{\zeta_0}.$$ And since the map $\Psi$
conserves the cross-ratio, we have  
$({\mathcal Z},\zeta_0+\zeta_1,\infty,\zeta_1)=(\cA,\cC,\cD,\cB)$.

Assume that $\cD=\cB$. Set $\zeta_1=\Psi(\cC)$. Then
$\zeta_1\neq\infty$ (since $\cC\neq\cD$), and $\Psi$
maps ${\mathcal A}$ to ${\mathcal Z}$, ${\mathcal B}$ to $\infty$, 
${\mathcal C}$ to $\zeta_1$ and ${\mathcal D}$ to $\infty$. Hence the
flux $\varphi$ of $Y$ through $E$ is equal to the flux of the Killing
field associated to the translation along the geodesic
$(\zeta_1,\infty)$ through 
$\Psi(E)$. This flux has been calculated in lemma
\ref{fluxcatends} : we have
$\varphi=\pi(1-\mu^2)$. And since $(\cA,\cC,\cD,\cB)=1$ in this
case, the result is still true. 

Assume that $\cC=\cB$. The flux with respect to the geodesic
$(\cC,\cD)$ is the opposite of the flux with respect to
$(\cD,\cC)$. Hence we have $\varphi=-\pi(1-\mu^2)$ according to what
has just been done. Consequently, since $(\cA,\cC,\cD,\cB)=0$ in this
case, the result is still true. 

We proceed in the same way for the flux of Killing fields associated to rotations. Then the expression of the flux polynomial follows from theorem \ref{fluxpolynomial}.
\end{proof}

\subsection{Flux for horospherical ends}

\begin{lemma} \label{fluxhoroends}
Let $E$ be a canonical horospherical end of asymptotic boundary $\infty$. Let
$\zeta_0$ and $\zeta_1$ be two complex numbers, with $\zeta_0\neq 0$.
Let $q_{-1}$ be the coefficient of the term of order $-1$ in the canonical
Hopf differential of the end.
Then the flux polynomial of $E$ is
$$\Pi_E(X)=-2\pi q_{-1}^2,$$
the flux of the Killing field associated to the translation along
the geodesic $(\zeta_0+\zeta_1,\zeta_1)$ through $E$ is 
$$-2\pi\re\left(\frac{q_{-1}^2}{\zeta_0}\right),$$
the flux of the Killing field associated to the rotation about
the geodesic $(\zeta_0+\zeta_1,\zeta_1)$ through $E$ is 
$$2\pi\im\left(\frac{q_{-1}^2}{\zeta_0}\right),$$
the flux of the Killing field associated to the translation along
the geodesic $(\zeta_1,\infty)$ through $E$ is zero,
and the flux of the Killing field associated to the rotation about
the geodesic $(\zeta_1,\infty)$ through $E$ is zero.
\end{lemma}

\begin{proof}
Using the canonical Bryant representation (\ref{hororepbryantcanon}), we compute that the coefficients of the flux polynomial are
$$\varphi_0=4\pi\Res(D\rmd C-C\rmd D)=-8\pi cg_1'(0),$$
$$\varphi_1=4\pi\Res(C\rmd B-D\rmd A)=0,$$
$$\varphi_2=4\pi\Res(B\rmd A-A\rmd B)=0.$$
Using equations (\ref{horof2prime}) and (\ref{horofg}), we obtain that $\varphi_0=-4\pi h'(0)$. Then we deduce from equations (\ref{horohprime}) and (\ref{horohopf}) that $\varphi_0=-2\pi q_{-1}^2$.
Applying theorems \ref{fluxtranslation}, \ref{fluxrotation} and \ref{fluxpolynomial}, we obtain the announced results.
\end{proof}

\begin{thm} \label{theoremfluxhorospherical}
Let $\cB\in\bar{\C}$ and $E$ be a
horospherical end of asymptotic boundary ${\mathcal B}$.

If $\cB\in\C$, then there exists a complex number $\kappa$ such that,
for all couples 
$(\cC,\cD)$ of distinct points in $\bar{\C}$, the flux of
the Killing field associated to the 
translation along the geodesic $({\mathcal C},{\mathcal D})$ through $E$ is
$$\varphi=-2\pi\re\left(\kappa
\frac{(\cC-\cB)(\cD-\cB)}{\cC-\cD}\right),$$
the flux of the Killing field associated to the 
rotation about the geodesic $({\mathcal C},{\mathcal D})$ through $E$ is
$$\varphi=2\pi\im\left(\kappa
\frac{(\cC-\cB)(\cD-\cB)}{\cC-\cD}\right),$$
and the flux polynomial of $E$ is 
$$\Pi_E(X)=-2\pi\kappa(X-\cB)^2.$$

If $\cB=\infty$, then there exists a complex number $\kappa$ such that,
for all couples 
$(\cC,\cD)$ of distinct points in $\bar{\C}$, the flux of
the Killing field $Y$ associated to the 
translation along the geodesic $({\mathcal C},{\mathcal D})$ through $E$ is
$$\varphi=-2\pi\re\left(\frac{\kappa}{\cC-\cD}\right),$$
the flux of the Killing field $Y$ associated to the 
rotation about the geodesic $({\mathcal C},{\mathcal D})$ through $E$ is
$$\varphi=2\pi\im\left(\frac{\kappa}{\cC-\cD}\right),$$
and the flux polynomial of $E$ is $$\Pi_E(X)=-2\pi\kappa.$$

The number $\kappa$ is called the flux coefficient of $E$. We have
$\kappa=0$ (or, equivalently, $\Pi_E(X)=0$) if and only if the Hopf differential $\omega\rmd g$ of the
end $E$ is holomorphic at zero, {\it i.e.} the degree $\mu$ of the
secondary Gauss map $g$ at zero is at least $3$.
\end{thm}

\begin{proof}
We first compute the flux $\varphi$ of the Killing field $Y$ associated to the
translation along the geodesic $({\mathcal C},{\mathcal D})$.

According to what has been done in section \ref{horoends}, there exists
a direct isometry $\Psi$ of $\h^3$ which maps $E$ to a canonical
horospherical end of asymptotic boundary $\infty$. We use the
notations of the beginning of section \ref{horocanrep}, with
$f_1'(0)=0$ (see remark \ref{derivzerohoro}).

Assume that neither $\cC$ nor $\cD$ is equal to $\cB$. Set
$\zeta_0=\Psi({\mathcal C})-\Psi({\mathcal D})$ and $\zeta_1=\Psi({\mathcal D})$.
Then $\zeta_0$ and $\zeta_1$ are different from $\infty$,
and $\Psi$ maps ${\mathcal B}$ to $\infty$,
${\mathcal C}$ to $\zeta_0+\zeta_1$ and ${\mathcal D}$ to $\zeta_1$. Hence the
flux $\varphi$ of $Y$ through $E$ is equal to the flux of the Killing
field associated to the 
translation about the geodesic $(\zeta_0+\zeta_1,\zeta_1)$ through
$\Psi(E)$. This flux has been calculated in lemma
\ref{fluxhoroends} : we have
$$\varphi=-2\pi\re\left(\frac{q_{-1}^2}{\zeta_0}\right).$$

We have $a=a_1f_1+a_2zf_2$ and $c=c_1f_1+c_2zf_2$. Hence we get
$a(0)=a_1$, $a'(0)=a_2$, $c(0)=c_1$ and $c'(0)=c_2$.

Using what has been done in the beginning of section \ref{horocanrep},
we compute that, if  $\cB\in\C$, then
$$\zeta_0=\frac{(a'(0)\cB-c'(0))^2(\cC-\cD)}{(\cC-\cB)(\cD-\cB)}$$ 
and if $\cB=\infty$, then $$\zeta_0=a'(0)^2(\cC-\cD).$$ 

We deal with the cases where $\cC$ or $\cD$ is equal to $\cB$ as for
theorem \ref{theoremfluxcatenoidal}, using lemma
\ref{fluxhoroends}.

We proceed in the same way for the flux of Killing fields associated to rotations. Then the expression of the flux polynomial follows from theorem \ref{fluxpolynomial}.

Moreover, the nullity of $\kappa$ is equivalent to the nullity of $q_{-1}$.
\end{proof}

\subsection{Flux for horospheres} 

\begin{thm} \label{transhorosphere}
Let $E$ be an
end which is part of a horosphere. Then the flux of the Killing field associated to the translation along any geodesic or to the rotation about any geodesic is zero, and the flux polynomial of $E$ is zero.
\end{thm}

\begin{proof}
Let $\Gamma$ be a generator of $\pi_1(E)$. Since a horosphere is
simply connected, $\Gamma$ is homotopic to zero in the
horosphere. Consequently, the fluxes are zero. Thus the flux polynomial is also zero.
\end{proof}

\section{Geometric applications}

\begin{defn}
Let $n$ be a positive integer. Let $\Sigma$ be a complete immersed Bryant
surface. We say that $\Sigma$ is a $n$-catenoidal surface if $\Sigma$ 
has exactly $n$ ends and each end is an embedded end of
finite total curvature.
\end{defn}

%\begin{thm} \label{thmpoly}
%Let $E$ be an embedded Bryant surface end of finite total
%curvature. Then there exists a unique polynomial $P_E(X,Y)\in\C[X,Y]$ such
%that, for all couples $(\cC,\cD)$ of distinct complex numbers, the flux
%of the Killing field associated to the translation along the geodesic
%$(\cC,\cD)$ through $E$ is $$\re\left(\frac{P_E(\cC,\cD)}{\cC-\cD}\right)$$
%and the flux of the Killing field associated to the rotation about
%the geodesic $(\cC,\cD)$ through $E$ is
%$$-\im\left(\frac{P_E(\cC,\cD)}{\cC-\cD}\right)$$
%
%This polynomial $P_E$ is symmetric.
%
%The polymomial $$\Pi_E(X)=P_E(X,X)$$ is called the flux polynomial of
%$E$. Its degree is at most $2$.
%
%If $E$ is a catenoidal end of growth $1-\mu$, of asymptotic boundary
%$\cB$ and of axis $(\cA,\cB)$, then
%$$\Pi_E(X)=2\pi(1-\mu^2)\frac{(X-\cA)(X-\cB)}{\cB-\cA}.$$
%(In the case where $\cA=\infty$ (respectively $\cB=\infty$), the above
%formula means $\Pi_E(X)=2\pi(1-\mu^2)(X-\cB)$ (respectively
%$\Pi_E(X)=-2\pi(1-\mu^2)(X-\cA)$).)
%
%If $E$ is a horospherical end of asymptotic boundary $\cB\in\C$ and of flux
%coefficient $\kappa$, then
%$$\Pi_E(X)=-2\pi\kappa(X-\cB)^2.$$
%
%If $E$ is a horospherical end of asymptotic boundary $\infty$ and of flux
%coefficient $\kappa$, then
%$$\Pi_E(X)=-2\pi\kappa.$$
%If $E$ is part of a horosphere, then $$\Pi_E(X)=0.$$
%\end{thm}
%
%\begin{proof}
%This is a reformulation of theorems \ref{theoremfluxcatenoidal},
%\ref{theoremfluxhorospherical}, \ref{transhorosphere},
%\ref{theoremfluxcatenoidalrot}, \ref{theoremfluxhorosphericalrot} and
% \ref{theoremfluxhororot}. 
%\end{proof}

\begin{prop} 
Let $\Sigma$ be a $n$-catenoidal surface. Then the sum of the fluxes
of any Killing field through its ends is zero.
\end{prop}

\begin{proof}
Let $W$ be a compact set in $\h^3$ such that $\Sigma\setminus W$ is
the disjoint union of the ends $E_j$ of $\Sigma$ and such that
$\partial W$ is a regular surface. Let $U_j$ be the part of 
$\partial W$ that is in the interior of $E_j$. Let $\Sigma'$ be the union of
$\Sigma\cap W$ and the $U_j$. We can calculate the flux of $E_j$ using
the curve $\partial U_j$ and the surface $U_j$. Since $\Sigma'$ is
homologous to $0$, the result follows from \cite{kkms}.
\end{proof}

\begin{cor} \label{poly}
Let $\Sigma$ be a $n$-catenoidal surface. Then the sum of the flux polynomials
of its ends is zero.
\end{cor}

\begin{prop} \label{deux}
Let $\Sigma$ be a $2$-catenoidal surface.
Assume that its ends $E_1$ and
$E_2$ are catenoidal ends of growths $1-\mu_1$ and $1-\mu_2$, of
asymptotic boundaries 
${\mathcal B}_1$ and ${\mathcal B}_2$,
and of axes $({\mathcal A}_1,{\mathcal B}_1)$ and
$({\mathcal A}_2,{\mathcal B}_2)$. Assume that $\cB_1\neq\cB_2$.
Then we have $\mu_1=\mu_2$, ${\mathcal A}_1={\mathcal B}_2$ and 
${\mathcal A}_2={\mathcal B}_1$ (that is to say, the two ends have the
same growth, the same axis, but two different asymptotic boundaries).
\end{prop}

\begin{proof}
Without loss of generality, we can assume that $\cA_1$, $\cA_2$, $\cB_1$ and
$\cB_2$ are different from $\infty$.

The sum of the flux polynomials of the two ends is zero. In
particular, these polynomials have the same roots. And since
$\cB_1\neq\cB_2$ we have $(\cA_1,\cB_1)=(\cB_2,\cA_2)$. Finally we
obtain $1-\mu_1^2=1-\mu_2^2=0$, {\it i.e} $\mu_1=\mu_2$.
\end{proof}

\begin{rem}
Levitt and Rosenberg (\cite{levitt}) have shown that if moreover
$\Sigma$ is properly embedded, then $\Sigma$ is a surface of
revolution, hence a catenoid cousin. It is essential that the end
should be properly
embedded : indeed Rossman and Sato (\cite{sato}) have
constructed a one-parameter family of genus one $2$-catenoidal
surfaces.
\end{rem}

\begin{rem}
The flux polynomial does not allow us to eliminate the case of a
$2$-catenoidal surface with two catenoidal ends having the same
asymptotic boundary. We do not know if such a surface exists. If it
exists, its ends have the same axis.
\end{rem}

\begin{prop} \label{trois}
Let $\Sigma$ be a $3$-catenoidal surface. Assume that its three ends are
catenoidal and that their asymptotic boundaries are distinct. Then,
given the growths, the axes of the three ends are uniquely determined,
they lie in the same plane and they are concurrent (possibly in the
asymptotic boundary of $\h^3$).
\end{prop}

\begin{proof} We use obvious notations.
Up to an isometry of $\h^3$, we can assume that $\cB_1=-1$, $\cB_2=0$ and
$\cB_3=1$. We set $\sigma_j=1-\mu_j^2$. Considering the coefficients
of the sum of the flux polynomials of the ends, we get
$$\frac{\sigma_1}{\cA_1+1}+\frac{\sigma_2}{\cA_2}
+\frac{\sigma_3}{\cA_3-1}=0,$$
$$\sigma_1\frac{\cA_1-1}{\cA_1+1}+\sigma_2+
\sigma_3\frac{\cA_3+1}{\cA_3-1}=0,$$
$$-\sigma_1\frac{\cA_1}{\cA_1+1}+\sigma_3\frac{\cA_3}{\cA_3-1}=0.$$ 

A computation gives
$$\cA_1=\frac{\sigma_1-\sigma_2+\sigma_3}{3\sigma_1+\sigma_2-\sigma_3},$$
$$\cA_2=\frac{\sigma_2}{\sigma_3-\sigma_1},$$
$$\cA_3=\frac{\sigma_1-\sigma_2+\sigma_3}{\sigma_1-\sigma_2-3\sigma_3}.$$

Consequently the points $\cA_j$ are uniquely determined. Moreover, all
the $\cA_j$ and $\cB_j$ are real. This means they lie in the same
plane. 

All the geodesics $(\cA_j,\cB_j)$ lie in the plane
$\{v=\im\zeta=0\}$. Assume that $\cA_1$, $\cA_2$ and $\cA_3$ are all
different from $\infty$. Then the equations of the geodesics $(\cA_1,-1)$,
$(\cA_2,0)$ and $(\cA_3,1)$ are respectively
$$u^2-u(\cA_1-1)+w^2-\cA_1=0,$$
$$u^2-u\cA_2+w^2=0,$$
$$u^2-u(\cA_3+1)+w^2+\cA_3=0.$$
Thus the abscissa of the intersection point of the first and the
second axes is $$u=\frac{\cA_1}{1-\cA_1+\cA_2},$$
and the abscissa of the intersection point of the second and the
third axes is $$u=\frac{\cA_3}{1-\cA_2+\cA_3}.$$
Hence the three axes are concurrent if and only if these two numbers
are equal. The expressions of the $\cA_j$ computed above show that
this is the case. 

We proceed in the same manner if exactly one of the $\cA_j$ is equal
to $\infty$. If two of the $\cA_j$ are equal
to $\infty$, then we deduce from the expressions of the $\cA_j$ that
the third one is also equal to $\infty$ ; in this case the axes are
concurrent at $\infty$.
\end{proof}

\begin{rem}
Levitt and Rosenberg (\cite{levitt}) have shown that if moreover
$\Sigma$ is properly embedded, then the plane $(\cB_1,\cB_2,\cB_3)$ is
a plane of symmetry of $\Sigma$ ; we can deduce from this that the
axes lie in this plane.
\end{rem}

There is an analogue of proposition \ref{trois} for minimal surfaces
in Euclidean space $\R^3$.

\begin{prop} \label{troisminimale}
Let $\Sigma$ be a minimal surface in $\R^3$. Assume that $\Sigma$ has
finite total curvature, three ends and that all the ends are
asymptotic to catenoids. Then the axes of the ends lie in the same
plane and they are either parallel or concurrent.
\end{prop}

\begin{proof}
Let $E_1$, $E_2$, $E_3$ be the ends of $\Sigma$, $F_j$ the flux of
$E_j$ and $T_j(P)$ its torque at the point $P$. We recall that the
axis of $E_j$ is the set of the points where $T_j$ is zero, and that
we have the formula $T_j(Q)=T_j(P)+\overrightarrow{PQ}\times
F_j$. Moreover, the two following ``balancing formulae'' hold :
$$F_1+F_2+F_3=0,$$
$$\forall P\in\R^3,T_1(P)+T_2(P)+T_3(P)=0.$$

We can assume that the three axes are not all identical (otherwise the 
result is clear). For each $j\in{1,2,3}$, let $P_j$ be a point of the
axis of $E_j$. Then the three axes are the straight lines $P_j+\R
F_j$. We can assume that $P_1$, $P_2$ and $P_3$ do not lie on the same 
straight line, since the axes are distinct.

We have $0=T_1(P_1)=-T_2(P_1)-T_3(P_1)=\overrightarrow{P_1P_2}\times F_2
+\overrightarrow{P_1P_3}\times F_3$. Hence we get
$$0=<T_1(P_1),\overrightarrow{P_1P_3}>=
<\overrightarrow{P_1P_2}\times F_2,\overrightarrow{P_1P_3}>=
\det(\overrightarrow{P_1P_2},F_2,\overrightarrow{P_1P_3}).$$ This
means that the axis of $E_2$ lies in the plane containing $P_1$, $P_2$ 
and $P_3$. We obtain the same result for $E_1$ and $E_3$. Hence the
three axes are coplanar.

If the axes of $E_1$ and $E_2$ are parallel, then the axis of $E_3$ is 
also parallel to them since $F_3=-F_1-F_2$.

If the axes of $E_1$ and $E_2$ are not parallel, then they meet at a
point $P_0$. Then we get $T_3(P_0)=-T_1(P_0)-T_2(P_0)=0$. So the three 
axes are concurrent at $P_0$.
\end{proof}

\begin{prop}
There is no $2$-catenoidal surface with one catenoidal end and one
horospherical end.
\end{prop}

\begin{proof}
Assume that such a surface exists. Without loss of generality, we can
assume that $\infty$ is not in the asymptotic boundary of the surface
and that $\infty$ is not an extremity of the axis of the
catenoidal end. Then the flux polynomials of its ends
have the same roots. This is impossible since the flux polynomial of a
catenoidal end has two simple roots and the flux polynomial of a
horospherical is either zero or has a double root.
\end{proof}

\begin{rem}
According to \cite{collin} (theorem 12), if a catenoidal surface is
properly embedded, then either it is a horosphere or all its ends are
catenoidal.
\end{rem}

\begin{ex}
In \cite{sousa}, de Sousa Neto has constructed Costa-type Bryant
surfaces. Let $\Sigma$ be such a surface. It is a $3$-catenoidal
surface of positive genus. It has two catenoidal ends $E_1$ and $E_2$,
which have the same asymptotic boundary $\cB$, and one horospherical end
$E_3$, whose asymptotic boundary $\cB_3$ is different from $\cB$. Let
$(\cA_1,\cB)$ and $(\cA_2,\cB)$ be the axes of $E_1$ and $E_2$, and $1-\mu_1$
and $1-\mu_2$ their respective growths. Since the sum of the flux
polynomial of the ends is zero, and since $\cB\neq\cB_3$, we obtain that
the flux coefficient of $E_3$ is zero, and then that $\cA_1=\cA_2$ and
$\mu_1^2+\mu_2^2=2$.
\end{ex}

\nocite{*}

\bibliographystyle{amsalpha}
\bibliography{flux}

\end{document}